\documentclass[12pt]{amsart}

\usepackage{amssymb}
\usepackage{a4wide}
\usepackage{graphicx}
\usepackage{tikz-cd}

\usepackage{subfigure}

\newcommand{\br}{{\rm br}}
\newcommand{\tr}{{\rm tr}}
\newcommand{\C}{\mathbb{C}}
\newcommand{\R}{\mathbb{R}}
\newcommand{\Q}{\mathbb{Q}}
\newcommand{\Z}{\mathbb{Z}}
\newcommand{\N}{\mathbb{N}}
\newcommand{\B}{\mathbb{B}}
\newcommand{\Sc}{\mathcal{S}}
\newcommand{\T}{\mathcal{T}}

\renewcommand{\k}{\mathbb{K}}
\renewcommand{\l}{\mathbb{L}}

\newcommand{\CHL}{\mathcal{C}}
\newcommand{\CP}{P_{\mathbb{C}}}
\renewcommand{\P}{\mathbb{P}}

\newcommand{\re}{\mathfrak{Re}}
\newcommand{\im}{\mathfrak{Im}}

\newtheorem{thm}{Theorem}[section]
\newtheorem{prop}{Proposition}[section]
\theoremstyle{remark}

\def\cqfd{\mbox{}\nolinebreak\hfill$\Box$\medbreak\par}
\newenvironment{pf}{\noindent\textbf{Proof:}}{\cqfd}

\usepackage{color}

\numberwithin{equation}{section}

\title{A new non-arithmetic lattice in $PU(3,1)$}

\author{Martin Deraux}

\date{Jul 1, 2019}

\begin{document}

\begin{abstract}
We study the arithmeticity of the Couwenberg-Heckman-Looijenga
lattices in $PU(n,1)$, and show that they contain a non-arithmetic
lattice in $PU(3,1)$ which is not commensurable to the non-arithmetic
Deligne-Mostow lattice in $PU(3,1)$. 
\end{abstract}

\maketitle

\section{Introduction}

Lattices in the isometry groups of most symmetric spaces of
non-compact type are arithmetic, due to celebrated superrigidity
results by Margulis~\cite{margulis} (symmetric spaces of higher rank),
Corlette~\cite{corlette}, and Gromov-Schoen~\cite{gromovschoen}
(quaternionic hyperbolic spaces, and the octonionic hyperbolic plane).

For small values of $n$, it is fairly easy to construct non-arithmetic
lattices in $SO(n,1)$ by using Coxeter polyhedra (a criterion due to
Vinberg gives a simple computational way to determine the
arithmeticity of these groups). For $n$ large enough, there are no
Coxeter polytopes in $H^n_{\mathbb{R}}$, but there are non-arithmetic
lattices in $SO(n,1)$ for arbitrary $n$ by a beautiful construction
due to Gromov and Piatetski-Shapiro~\cite{gps}. Their construction
produces infinitely many commensurability classes of non-arithmetic
lattices in any dimension. Note that the general structure of lattices
in $SO(n,1)$ remains mysterious.

The situation is even more mysterious for lattices in $PU(n,1)$,
$n\geq 2$, which is (up to finite index) the isometry group of complex
hyperbolic space $H^n_{\mathbb{C}}$. Here there is currently no
analogue of the Gromov-Piatetski-Shapiro construction (there exist no
real totally geodesic hypersurfaces in $H^n_{\mathbb{C}}$, so there is
no reasonable gluing interface to construct hybrids). In fact, only
finitely many commensurability classes of non-arithmetic lattices in
$PU(n,1)$ are known, only for very low values of $n$.

The first examples in $PU(2,1)$ were due to
Mostow~\cite{mostowpacific}, and his construction was soon
generalized to produce several more examples in $PU(2,1)$, and a
single one in $PU(3,1)$ see~\cite{delignemostow}. For some decades,
the Deligne-Mostow examples were the only known examples, even though
some alternative constructions were given, see~\cite{thurstonshapes}
for instance. To this day, it is still unknown whether there
exist non-arithmetic lattices in $PU(n,1)$ for any $n>3$.

A slightly different construction was given by Hirzebruch
(see~\cite{bhh}), based on the equality case in the Miyaoka-Yau
inequality, i.e. an orbifold version of the fact that a compact
complex surface $X$ of general type with $c_1^2(X)=3c_2(X)$ is covered
by the ball. Given such an $X$, the existence of a lattice $\Gamma$ in
$PU(2,1)$ such that $X=\Gamma\setminus \B^2$ is guaranteed, but it is
not obvious how to describe the lattice explicitly (the existence of a
K\"ahler-Einstein metric is obtained by showing existence of a
solution to a Monge-Amp\`ere equation).

In fact, the arithmetic structure of the Hirzebruch examples seems not
to have been worked out anywhere in the literature, apart from a small
number of examples where coincidences with some arithmetic groups were
found (see the work of
Holzapfel~\cite{holzapfel},~\cite{holzapfelbook}, and also the more
recent~\cite{derauxklein}).

The Deligne-Mostow construction and the Barthel-Hirzebruch-H\"ofer
construction were given a common generalization by Couwenberg, Heckman
and Looijenga~\cite{chl}, but their work barely brushes the discussion
of arithmeticity (they mention that the examples derived from real
Coxeter groups are arithmetic, without any details). It was recently
observed~\cite{derauxklein} that some of the non-arithmetic lattices
in $PU(2,1)$ produced by the author, Parker and Paupert~\cite{dpp2}
were in fact conjugate to some specific Couwenberg-Heckman-Looijenga
lattices.

The main goal of the present paper is to give a systematic study of
the arithmeticity of the Couwenberg-Heckman-Looijenga lattices.  We
write $\CHL(G,p_1,\dots,p_k)$ for the CHL lattice derived from the
Shephard-Todd group $G$, generated by complex reflections of angle
$2\pi/p_j$. As mentioned in~\cite{chl}, when $G$ is the Weyl group of
type $A_n$ or $B_n$, the lattices of the form $\CHL(G,p)$ are all
commensurable to Deligne-Mostow lattices (the Deligne-Mostow
construction gives lattices in $PU(n,1)$ only for $n\leqslant
9$). Note also that the Shephard-Todd group $G_{32}$ is obtained from
the Couwenberg-Heckman-Looijenga construction starting with the group
$W(A_4)$ so, just as in~\cite{chl}, we omit $G_{32}$ from our lists
(it would again produce lattices commensurable with Deligne-Mostow
lattices). The imprimitive Shephard-Todd groups $G(m,p,n)$ can also be
obtained from classical groups by the CHL construction, so the
corresponding ball quotients are also commensurable to Deligne-Mostow
ball quotients.

We refer to (primitive) Shephard-Todd groups not of type $A_n$ or
$B_n$ as \emph{exceptional complex reflection groups}, and we refer to
the corresponding mirror arrangements as \emph{exceptional
  arrangements}. Via the CHL construction, the exceptional
arrangements produce lattices in $PU(n,1)$ only for $n\leqslant 7$.

The CHL lattices in $PU(2,1)$ were already mentioned
in~\cite{derauxklein} and~\cite{derauxabelian}. It turns out that, for
$n\geq 3$, all the (non Deligne-Mostow) CHL lattices in $PU(n,1)$ are
arithmetic except for one.
\begin{thm}
  Let $\Gamma$ be a CHL lattice derived from an exceptional finite complex
  reflection group $G$ acting irreducibly on $V=\C^{n+1}$, $n \geqslant
  3$. Then $\Gamma$ is arithmetic, unless $\Gamma=\CHL(G_{29},3)$.
\end{thm}

More precisely, we state the following.
\begin{thm}\label{thm:new3d}
  The lattice $\CHL(G_{29},3)$ is a non-arithmetic, non-cocompact
  lattice, with adjoint trace field $\Q(\sqrt{3})$. It is not
  commensurable to any Deligne-Mostow lattice.
\end{thm}
Recall that the Deligne-Mostow list of lattices contains only one
non-arithmetic lattice in $PU(n,1)$ with $n\geq 3$, namely the lattice
$\Gamma_{\mu}$ for $\mu=(3,3,3,3,5,7)/12$; so the main additional
content of Theorem~\ref{thm:new3d} is the claim that $\CHL(G_{29},3)$
is not commensurable to that specific $\Gamma_\mu$.

Given the commensurability analysis in~\cite{thealgo}, putting
together all known non-arithmetic lattices in $PU(n,1)$, we see that
there are currently 22 known commensurability classes in $PU(2,1)$ and
2 commensurability classes in $PU(3,1)$.

The basic tool for proving these results is the knowledge of explicit
presentations of the braid groups associated to the Shephard-Todd
groups (see the conjectural statements in~\cite{bessismichel}, later
proved in~\cite{bessisannals}).  Using braid relations between the
generators, we study the irreducible representations of the
corresponding braid groups that send the generators to complex
reflections of the appropriate angle (the values of the angle for the
discrete holonomy groups in Couwenberg-Heckman-Looijenga have been
tabulated, see section~8 of~\cite{chl}).

It turns out there are finitely many such representations, and the
finite number is usually very small. Basic geometric considerations
(using cocompactness or discreteness arugments) allow us to single out
(a group conjugate to) the Couwenberg-Heckman-Looijenga holonomy
group.  Along the way, we find explicit matrices for generators for
the holonomy groups, which may be of independent interest (but these
were not given in~\cite{chl}).

Each holonomy group preserves an explicit Hermitian form (by
irreducibility, such an invariant Hermitian form is unique up to
scaling). The strategy for determining arithmeticity is then to
\begin{enumerate}
\item Find coordinates such that the Hermitian form has entries in a
  number field;
\item Check that the above number field is as small as possible;
\item Find coordinates where the matrices of the generators are
  actually algebraic integers.
\end{enumerate}
It is known that~(1) can always be achieved, because of Calabi-Weil
local rigidity of lattices, see chapter~VI of~\cite{raghunathan}. In
general, it is not easy to make that result effective, but it turns
out to be fairly easy in the cases we consider in the paper.

Step~(2) follows from the determination of the adjoint trace field,
i.e. the field generated by traces in the adjoint representation,
which is a well known commensurability invariant for lattices (in fact
for Zariski dense groups).

It is not known whether step~(3) can always be achieved, even though
it is strongly believed to be the case for every lattice in $PU(n,1)$
(for cocompact lattices, it follows from very recent work of Esnault
and Groechenig~\cite{esnault}). Recall that there are so-called
quasi-arithmetic lattices in $SO(n,1)$ for every $n$, i.e. lattices
where arithmeticity fails only by failure of integrality
(see~\cite{belothomson} and~\cite{thomson}).

We will go through steps~(1) through~(3) by explicit case by case
computation. In fact, we follow a suggestion of the referee and
combine steps~(2) and~(3).

Some parts of the paper require delicate arguments. One is the proof
that $\CHL(G_{29},3)$ is not commensurable to the Deligne-Mostow
non-arithmetic lattice in $PU(3,1)$. Indeed, the two groups have the
same rough commensurability invariants (cocompactness, adjoint trace
field and non-arithmeticity index, as defined in section 6.2
of~\cite{thealgo}). We work out an explicit description of the cusps
of these two lattices, and show that the cusps themselves are not
commensurable. Another delicate part is the determination of the
reflection representations for the braid group associated to the
Shephard-Todd group $G_{31}$. This group is not well-generated, in the
sense that it is not generated by the right number of reflections for
the ambient dimension.

\noindent \textbf{Acknowledgements:} It is a pleasure to thank
St\'ephane Druel, John Parker, Erwan Rousseau and Domingo Toledo for
their enthusiasm about this project. I am indebted to Gert Heckman and
Eduard Looijenga for explaining certain points in~\cite{chl}. I am
also very greatful to the referee, who suggested several significant
improvements of the manuscript.

\section{Basic facts about complex reflections.}\label{sec:basics}

We start with a complex vector space $V$ equipped with a
non-degenerate Hermitian inner product, which we denote by $\langle
\cdot,\cdot\rangle$ (we take this to be linear on the first factor,
and antilinear on the second factor). A \emph{complex reflection} is
a linear transformation of the form $R_{v,z}$ where
\begin{equation}
  R_{v,z}(x) = x+(z-1)\frac{\langle x,v\rangle}{\langle v, v\rangle}v,\label{eq:refl}
\end{equation}
for some vector $v\in V$ with $\langle v,v\rangle\neq 0$, and some
$z\in\C$ with $|z|=1$. It is easy to see that such a
transformation preserves the Hermitian inner product.

Note that scaling the vector $v$ does not change the above
transformation, so $v$ is not uniquely determined by the
transformation. The reflection fixes pointwise the complex-linear
subspace $v^\perp=\{w\in\C^{n,1}:\langle w,v\rangle=0\}$, called its
\emph{mirror}, and it acts on $\C v$ by multiplication by $z$. We
will call such a vector $v$ the \emph{polar vector} to the mirror
(this is only well-defined up to scaling). The complex number $z$
is called the \emph{multiplier}, and its argument is called the
\emph{angle} of the complex reflection.  In this paper, we will only
consider reflections of finite order, i.e. $z$ will actually be a
root of unity.

When the Hermitian inner product is positive definite, the condition
$\langle v,v\rangle\neq 0$ is of course equivalent to $v\neq 0$, and
the definition of a complex reflection agrees with the one
in~\cite{shephardtodd}. We will also use the same definition for
hyperbolic Hermitian inner products, i.e. those of signature $(n,1)$,
which are related to complex hyperbolic space $H^n_{\mathbb{C}}$ (for
basic information about complex hyperbolic space, see~\cite{goldman}
for instance). As a set, $H^n_{\mathbb{C}}$ is the set of complex
lines spanned by vectors $v$ with $\langle v,v\rangle<0$, and the
metric is built in such a way that the linear isometries of the
Hermitian inner product induce isometries of $H^n_{\mathbb{C}}\subset
\mathbb{P}(V)$ (in fact, the corresponding group $PU(n,1)$ has index
two in the full isometry group, the latter being obtained by adjoining
any antiholomorphic isometry).

In the hyperbolic case, we require moreover that $v$
in equation~\eqref{eq:refl} be a positive vector, i.e. $\langle
v,v\rangle>0$. In that case, the restriction of the Hermitian form to
$v^\perp$ has signature $(n-1,1)$, and the set of negative vectors in
$v^\perp$ projects down to a totally geodesic copy of
$H^{n-1}_{\mathbb{C}}$.  Moreover since we are free to scale $v$, we
can (and will often) assume $\langle v,v \rangle=1$.

For $k\in\N^*$, two group elements $a$ and $b$ are said to satisfy a
braid relation of length $k$ if
$$
 (ab)^{k/2}=(ba)^{k/2}.
$$ In that case, we write $\br_k(a,b)$. When $k$ is odd, the notation
$(ab)^{k/2}$ stands for an alternating product $a\cdot b\cdot a \cdots
b\cdot a$ with $k$ factors.  Note that when $\br_k(a,b)$ holds,
$\br_{nk}(a,b)$ also holds for every $n\geq 1$. The smallest $k$ such
that $\br_k(a,b)$ holds is called the braid length of $a$ and $b$,
which we denote by $\br(a,b)$.

It is often convenient to describe reflection groups by a (complex)
Coxeter diagram. The diagram is attached to a generating set of
reflections (with pairwise distinct mirrors). It has one vertex for
each generating complex reflection, and vertices are represented by a
circled integer, to indicate the order of the corresponding complex
reflection (more precisely, a node with a circled $p$ stands for a
complex reflection with multiplier $e^{2\pi i/p}$).

A pair of vertices is joined by an edge labelled $k$ if the
corresponding reflections satisfy a braid relation of length $k \geq
3$.  Braid relations of length 3 are called standard braid relations,
and the corresponding edge in the graph is drawn, but the label 3 is
usually omitted. Braid relations of length 4 are often drawn by
doubling the corresponding edge in the graph (and then omitting the
label 4).

In order to get abstract presentations for Shephard-Todd groups, we
often need to include extra relations, which are indicated with extra
decorations of the diagram, see the Appendix~2 in~\cite{brmaro}. For
instance, in the diagram~(d) in Figure~\ref{fig:cox}, the vertical
double bar indicates a braid relation $\br_4(R_3,R_2R_4)$.

\begin{figure}
  \includegraphics[width=5cm]{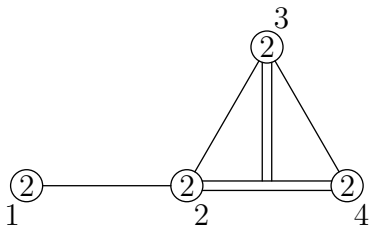}
  \caption{Coxeter diagram for $G_{29}$}\label{fig:cox-g29}
\end{figure}

The Coxeter diagrams for the Shephard-Todd groups that we will need in
this paper can be deduced from the graphs given in Table~\ref{fig:cox}
by replacing the nodes by circled 2's. For example, the Coxeter diagram
for $G_{29}$ is the one in Figure~\ref{fig:cox-g29}.
We have numbered the nodes, and write $R_j$ for the reflection
corresponding to the node labelled $j$ ($j=1,2,3,4$). The
corresponding presentation is the one given in
equation~\eqref{eq:braidgroup-g29} (p.~\pageref{eq:braidgroup-g29}),
with extra relations $R_j^2=Id$ (in fact only one of these four
relations suffices, since $\br_3(R_j,R_{j+1})$ implies that $R_j$ and
$R_{j+1}$ are conjugate).

We only list exceptional Shephard-Todd groups in $U(n+1)$ ($n\geqslant
3$) with a generating set consisting of reflections of order 2, since
the ones with higher order generators do not produce any more lattices
(see~\cite{chl}). 

We briefly sketch the general strategy for writing matrices for the
Hermitian forms and complex reflections generating CHL lattices. For
well-generated groups, i.e. reflection groups in $U(n+1)$ generated by
$n+1$ reflections, we will use coordinates given by the basis obtained
by choosing polar vectors for these $n+1$ reflections; in other words,
the polar vectors are simply given by the standard basis vectors
$e_1,\dots,e_{n+1}$ of $\C^{n+1}$. The reflections are uniquely
determined by the matrix of the Hermitian form in that basis, given by
$H_{j,k}=\langle e_k,e_j\rangle$.

We may assume that $\langle e_j,e_j\rangle=1$, so we will take $H$ to
have ones on the diagonal. For $j\neq k$, the braid length
$\br(R_j,R_k)$ determines $|\langle e_j,e_k\rangle|$ (see section~2.2
of~\cite{mostowpacific}). We can choose the argument of $\langle
e_j,e_k\rangle$ freely, by replacing $e_k$ by $\lambda e_k$ for some
$|\lambda|=1$, and the goal will be to get the corresponding
reflections to have algebraic integer entries (in the discussion
below, we will call a choice of argument \emph{reasonable} if the
corresponding matrices have algebraic integer entries).

For triples $j,k,l$ of pairwise distinct indices, we can freely choose
the arguments of $\langle e_j,e_k\rangle$ and $\langle
e_k,e_l\rangle$, but then the third one $\langle e_k,e_l\rangle$,
cannot usually be chosen arbitrarly, unless one of the inner product
is zero (for more details on this, see section~2.3
of~\cite{mostowpacific}).

An important special case is the one where the Coxeter diagram is a
tree (see (a), (b), (c), (e), (i), (j), (k) in
Figure~\ref{fig:cox}). In that case, the arguments of the inner
products $\langle e_j,e_k\rangle$ can be chosen freely (for instance
we may assume all inner products are real, but this is not always a
``reasonable'' choice, in the above sense).

In other well-generated cases (see (d), (g), (h) in
Figure~\ref{fig:cox}), we consider triangles in the diagram. If a
triangle corresponds to polar vectors $e_j,e_k,e_l$, we fix a
reasonable choice of $\langle e_j,e_k\rangle$ and $\langle
e_k,e_l\rangle$, and use the presentation of the braid group to find
admissible values of $\langle e_j,e_l\rangle$ (in well-generated
cases, there are only finitely many admissible values).

If the corresponding matrices have algebraic integer entries in the
correct number field, we are done. Otherwise, we use ad-hoc changes of
coordinates that produce algebraic integer entries (see
section~\ref{sec:g28} for the group $G_{28}$).

The only group we need to consider that is not well generated is
$G_{31}\subset U(4)$. In fact that group is generated by $5=4+1$
reflections. In that case, we will use the same parametrization as
above using only 4 of the 5 reflections, and then use braid relations
in the Bessis-Michel presentation to determine the polar vector to the
mirror of the fifth reflection. For more details on this, see
section~\ref{sec:g31}).

\section{The Couwenberg-Heckman-Looijenga lattices}

In~\cite{chl}, Couwenberg-Heckman-Looijenga give a general
construction of affine structures on the complement of hyperplane
arrangements in projective space, parametrized by angles related to
the holonomy around the hyperplanes in the arrangement. They also give
necessary and sufficient conditions for the completion of that
structure to be an orbifold (i.e. the holonomy is discrete, and the
completion is a quotient of the appropriate complex space form).

It is unclear how often these conditions are satisfied, but there is a
somewhat large list of examples associated to finite unitary groups
generated by complex reflections (these were classified by Shephard
and Todd~\cite{shephardtodd}). That list contains a lot of the
previously known examples of lattices in $PU(n,1)$ generated by
complex reflections, namely the Deligne-Mostow
lattices~\cite{delignemostow}, as well as the ones constructed by
Barthel, Hirzebruch and H\"ofer~\cite{bhh}. Note that some examples
in~\cite{thealgo} are still not covered by the
Couwenberg-Heckman-Looijenga construction, see~\cite{derauxabelian}.

The Couwenberg-Heckman-Looijenga lattices are decribed by giving:
\begin{itemize}
  \item an irreducible Shephard-Todd group $G$;
  \item a positive integer $p_j\geqslant 2$, $j=1,\dots,k$ for each
    of the $k$ orbits of mirrors of complex reflections in $G$.
\end{itemize}
We denote by $\CHL(G,p_1,\dots,p_k)$ the corresponding group. In this
paper, we only consider the exceptional Shephard-Todd groups, since
the other ones are covered by Deligne-Mostow theory.

It turns out (exceptional) Shephard-Todd groups have at most two
orbits of mirrors, so we only take $k\leq 2$. In fact, there is a
single orbit of mirrors (i.e. $k=1$) for all but one group, namely
$G=G_{28}$, which is isomorphic to the Coxeter group $F_4$.

The CHL structures are obtained as structures on the complement of the
union of the mirrors of reflections in $G\subset U(n+1)$; we follow
the notation in~\cite{chl} and write $V=\C^{n+1}$, $\mathcal{H}$ for
the union of mirrors of reflections in $G$, and $V^0=V\setminus
\mathcal{H}$.  In particular, by construction, the holonomy group is a
quotient of $\pi_1(\P(V^0/G))$, which is often called a \emph{braid
  group}.

It is reasonably easy (especially using modern computer technology,
and more so in low dimensions) to write down explicit group
presentations in terms of generators and relations for the
Shephard-Todd groups. This was done by Coxeter, see~\cite{coxeterbook}
and also Appendix 2 in~\cite{brmaro} for a convenient list.

This gives presentations for some quotients of the braid group
$\pi_1(\P(V^0/G))$, namely the orbifold fundamental group of the
quotient $\P(V/G)$, but it is not completely obvious how to deduce a
presentation for $\pi_1(\P(V^0/G))$. Roughly speaking, one would like
to cancel the relations expressing the order of reflections, and keep
the braid relations, but this is of course not
well-defined. Presentations for $\pi_1(\P(V^0/G))$ were proposed by
Bessis and Michel in~\cite{bessismichel}, and their conjectural
statements have later been proved in~\cite{bessisannals}.

Note that the Bessis-Michel presentations are given in such a way that
the generators correspond to suitably chosen simple loops around
hyperplanes in the arrangement. It follows that the
Couwenberg-Heckman-Looijenga holonomy groups $\CHL(G,p)$
(resp. $\CHL(G,p_1,p_2)$) are homomorphic images of the braid group
$\pi_1(\P(V^0/G))$, such that the corresponding homomorphism maps the
Bessis-Michel generators to complex reflections of angle $2\pi/p$
(resp. $2\pi/p_1$ and $2\pi/p_2$). One can in fact obtain an explicit
presentation for the lattices in terms of these generators (see
Theorem~7.1 in~\cite{chl} or section 4 in~\cite{chl3d_2}).

For arrangements of type $A_n$ or $B_n$, the corresponding lattices
are commensurable to Deligne-Mostow lattices, and the list is a bit
too long to be reproduced here, see p.157-159 of~\cite{chl}.
%, and also
%section~\ref{sec:dm} in this paper for the 3-dimensional case.
The other CHL lattices (in $PU(n,1)$ with $n\geqslant 3$) are listed in
Tables~\ref{tab:2d} and~\ref{tab:list} in appendix~\ref{sec:rough}.
%(bold-face means cocompact, red means
%non-arithmetic).

\section{Arithmeticity} \label{sec:criterion}

We will use the following arithmeticity criterion, which is proved
in~\cite{mostowpacific} (see also~\cite{delignemostow}). We refer to
it as the Vinberg/Mostow arithmeticity criterion. In what follows,
$\tr Ad \Gamma=\Q(\{\tr Ad\gamma:\gamma\in\Gamma\})$ is the field
generated by traces of elements of $\Gamma$ in the adjoint
representation.
\begin{thm}\label{thm:criterion}
  Let $H$ be a Hermitian form of signature $(n,1)$, defined over a CM
  field $\l\supset\k$. Let $\Gamma$ be a lattice
  in $SU(H,\mathcal{O}_\l)$, such that $\tr Ad
  \Gamma=\k$. Then $\Gamma$ is arithmetic if and only if $H^\sigma$ is
  definite for every $\sigma\in Gal(\l)$ acting
  non-trivially on $\k$.
\end{thm}
Recall that a CM field is a purely imaginary quadratic extension of a
totally real number field, we denote by $\k$ the totally real field
and by $\l$ the imaginary quadratic extension. As usual,
$\mathcal{O}_\l$ denotes the ring of algebraic integers. Note that not
every lattice is commensurable to a lattice as in the above statement,
which are sometimes called lattices \emph{of simplest type} (for the
general case, one needs to consider division algebras over a CM field).

Because of the fact that the adjoint representation of a unitary
representation $\rho$ is isomorphic to the tensor product
$\rho\otimes\overline{\rho}$, we have
$$
  \tr Ad \gamma=|\tr \gamma|^2
$$
for all $\gamma\in\Gamma$, which we will repeatedly use in the sequel.

\section{Explicit generators and arithmeticity} \label{sec:arithmeticity}

The goal of this section is to give explicit matrix generators for the
CHL lattices, as well as explicit Hermitian forms, and to use these to
apply the arithmeticity criterion stated in
section~\ref{sec:criterion}. We only work on lattices derived from
exceptional complex reflection groups acting on $\C^{n+1}$ with $n\geq
3$, which give an action on $\P^{n}$ with $n\geq 3$. Indeed,
non-exceptional ones yield Deligne-Mostow groups (explicit matrices
can easily be deduced from~\cite{delignemostow}, see
also~\cite{terada}); 2-dimensional examples turn out to be
commensurable to groups that have been studied elsewhere
(see~\cite{thealgo}, for instance).

We go through a somewhat painful case by case analysis in
sections~\ref{sec:g28} through~\ref{sec:g37}. The groups in
sections~\ref{sec:g28} through~\ref{sec:g31} give lattices in
dimension 3, the next ones in slightly higher dimension (the list of
groups and respective dimensions is given in
appendix~\ref{sec:rough}).

Note that, just as in~\cite{chl}, we do not include the group $G_{32}$
in the list, since it can be seen as a group derived from the $A_4$
arrangement, hence the corresponding lattices already appear in the
Deligne-Mostow list (see~\cite{chl3d_2} for more details).

\subsection{Lattices derived from $G_{28}$}\label{sec:g28}

Recall that $G_{28}$ has two orbits of mirrors of reflections
(see~\cite{chl} for instance), hence the corresponding CHL lattices
depend on two integer parameters. We denote the corresponding groups
by $\CHL(G_{28},p,q)$.

We call $r_1,\dots,r_4$ generators of $G_{28}$, numbered according to
the numbering of the nodes in Figure~\ref{fig:cox}(c)
(page~\pageref{fig:cox}). The orbits of mirrors can be checked
to be represented by the mirrors of $r_1$ and $r_4$.
%
%Indeed, $r_1$ and $r_2$ are conjugate since $\br_3(r_1,r_2)$
%holds. Similarly the mirror$r_3$ and $r_4$ are in the same orbit. It
%is not completely obvious, but it is a standard fact that $r_1$ and
%$r_4$ are not conjugate in the group (this can of course easily be
%checked with a computer algebra program, say GAP).

Couwenberg, Heckman and Looijenga show that there exist
representations of $G_{28}$ into $PU(3,1)$, with lattice image,
mapping $r_1,r_2$ to complex reflections with multiplier $e^{2\pi
  i/p}$ and $r_3,r_4$ to complex reflections of multiplier $e^{2\pi
  i/q}$, for $(p,q)$ given by $(2,q)$, $q=4,5,6,8,12$, $(3,q)$ for
$q=3,4,6,12$, $(4,4)$ and $(6,6)$.

For a generic value of $p,q$, we set up the Hermitian form as 
\begin{equation}\label{eq:generic_g28}
\left(\begin{matrix}
  1 & \alpha & 0 & 0\\
  {\bar\alpha} & 1 & \beta & 0\\
  0 & {\bar\beta} & 1 & \gamma\\
  0 & 0 & {\bar\gamma} & 1
\end{matrix}\right),
\end{equation}
and the generators are given by
$R_{e_1,z},R_{e_2,z},R_{e_3,w},R_{e_4,w}$, where $z=e^{2\pi i/p}$,
$w=e^{2\pi i/q}$ and the $e_j$, $j=1,2,3,4$ are the standard basis
vectors of $\C^4$. The corresponding matrices are given in
equation~\eqref{eq:firsttry}.
\begin{eqnarray}\label{eq:firsttry}
R_1=\left(\begin{matrix}
  z & \alpha(z-1) & 0 & 0\\
  0 & 1 & 0 & 0\\
  0 & 0 & 1 & 0\\
  0 & 0 & 0 & 1
\end{matrix}\right),
R_2=\left(\begin{matrix}
  1 & 0 & 0 & 0\\
  {\bar\alpha}(z-1) & z & \beta(z-1) & 0\\
  0 & 0 & 1 & 0\\
  0 & 0 & 0 & 1
\end{matrix}\right),\\
R_3=\left(\begin{matrix}
  1 & 0 & 0 & 0\\
  0 & 1 & 0 & 0\\
  0 & {\bar\beta}(w-1) & w & \gamma(w-1)\\
  0 & 0 & 0 & 1
\end{matrix}\right),
R_4=\left(\begin{matrix}
  1 & 0 & 0 & 0\\
  0 & 1 & 0 & 0\\
  0 & 0 & 1 & 0\\
  0 & 0 & {\bar\gamma}(w-1) & w
\end{matrix}\right).
\end{eqnarray}

The braid relation $\br_3(R_1,R_2)$ is equivalent to
$|\alpha|=\frac{1}{|z-1|}=\frac{1}{2\sin\frac{\pi}{p}}$
(see~\cite{mostowpacific} for instance). Similarly $\br_3(R_3,R_4)$ is
equivalent to $|\gamma|=\frac{1}{|w-1|}$.

One checks that $\br_4(R_2,R_3)$ is equivalent to $\beta=0$ or
$$
\beta^2=\frac{z+w}{z+w-1-zw}=\frac{\cos(\frac{\pi}{p}-\frac{\pi}{q})}{2\sin\frac{\pi}{p}\sin\frac{\pi}{q}}.
$$ 
We first rule out the case $\beta=0$.
\begin{prop}\label{prop:donotcommute}
  In the CHL lattice $\CHL(G_{28},p,q)$, $R_3$ and $R_4$ do not commute.
\end{prop}
\begin{pf}
  Let $S_2$ and $S_3$ denote reflections in $G_{28}$ acting on $\CP^3$
  that correspond to $R_2$ and $R_3$ in $\CHL(G_{28},p,q)$. By
  ``corresponding'', we mean that $R_j$ and $S_j$ are images of the
  same element $r_j$ in the braid group $\pi_1(P(V^0/G_{28}))$.

 Then $S_2$ and $S_3$ generate a group of order 8, isomorphic to the
 imprimitive Shephard-Todd group $G(4,4,2)$, and the arrangement has 4
 planes intersecting along the mirror intersection $L_2\cap L_3$,
 namely the mirrors of $S_2$, $S_3$, $S_2S_3S_2$, $S_3S_2S_3$.

 The branch locus of the quotient map $\C^2\rightarrow \C^2/G(4,4,2)$
 has local analytic structure $z_1^4=z_2^2$, which gives two tangent
 components (this can be seen by computing the invariant polynomial
 ring, see~\cite{bannai} for instance). This gives the structure of
 the quotient $\CP^3/G_{28}$, near a generic point of the intersection
 of the mirrors of $S_2$ and $S_3$, the branch locus is given locally
 analytically by the same equation $z_1^4=z_2^2$ (but there is a third
 variable, say $z_3$).

 If $R_2$ and $R_3$ were to commute, their mirrors would have to be
 orthogonal (they cannot coincide, otherwise the monodromy group would
 not act irreducibly on $\C^4$), and the branch locus in the quotient
 would have two smooth transverse components, which is a
 contradiction.
\end{pf}

We take $\alpha=\frac{1}{z-1}$ and $\gamma=\frac{1}{w-1}$ (this is a
natural choice, given the shape of the matrices for $R_1$ and $R_3$ in
equation~\eqref{eq:firsttry}), $\beta=\sqrt{\frac{z+w}{z+w-1-zw}}$,
and we take as the basis for $\C^4$ the vectors
$e_1,e_2,Je_2,J^{-1}e_2$ where $J=R_2R_3R_4$ (where the $e_j$ are
simply the standard basis vectors of $\C^4$).  These vectors do indeed
form a basis as long as $w+z\neq 0$, which will be the case for all
relevant pairs $(p,q)$.

We write $Q$ for the corresponding matrix
$$
Q = \left(\begin{matrix}
  1 & 0 & 0                      & 0\\
  0 & 1 & -w                     & \bar z\\
  0 & 0 & \frac{z+w}{\beta(1-z)} & \frac{z+w}{\beta zw(z-1)}\\
  0 & 0 & 0                      & \frac{z+w}{\beta zw(z-1)}
\end{matrix}\right),
$$
and get the matrices $\widetilde{R}_j=Q^{-1}R_jQ$ to have entries in $\Z[z,w]$, namely
\begin{eqnarray}\label{eq:secondtry}
&\widetilde{R}_1=
  \left(\begin{matrix}
    z & 1 & -w & \bar z\\
    0 & 1 & 0  & 0\\
    0 & 0 & 1  & 0\\
    0 & 0 & 0  & 1
  \end{matrix}\right),\quad 
\widetilde{R}_2=
  \left(\begin{matrix}
    1  & 0 & 0       & 0\\
    -z & z & -z(w+1) & 1+\bar w\\
    0  & 0 &  1      & 0\\
    0  & 0 &  0      & 1
  \end{matrix}\right)\\
&\widetilde{R}_3=
  \left(\begin{matrix}
    1 & 0   & 0  & 0\\
    0 & 1+w & -w & 0\\
    0 &  1  & 0  & 0\\
    0 &  0  & 0  & 1
  \end{matrix}\right),\quad
\widetilde{R}_4=
  \left(\begin{matrix}
    1 & 0 & 0    & 0\\
    0 & 1 & 0    & 0\\
    0 & 0 & 1+w  & -\bar z \bar w\\
    0 & 0 & zw^2 & 0
  \end{matrix}\right),
\end{eqnarray}
which preserve a Hermitian form defined over $\Q(z,w)$, in fact
\begin{equation}\label{eq:form-g28}
Q^*HQ=
\left(
\begin{matrix}
  1                          & \frac{1}{z-1}                     & \frac{w}{1-z}          & \frac{1}{z(z-1)}\\
  \frac{1}{\bar z-1}         &      1                            &  \frac{z(w+1)}{1-z}     & \frac{1+\bar w}{z-1}\\
  \frac{\bar w}{1-\bar z}    & \frac{\bar z(\bar w+1)}{1-\bar z}  &         1              & \frac{\bar w-\bar w^2}{1-z}\\
  \frac{1}{\bar z(\bar z-1)} & \frac{1+w}{\bar z-1}              & \frac{w-w^2}{1-\bar z} &        1
\end{matrix}
\right)
\end{equation}

\begin{prop} \label{prop:pre}
  The adjoint trace field of $\CHL(G_{28},p,q)$ is given by
  $\Q(\cos\frac{2\pi}{l})$, where $l$ is the least common multiple of
  $p$ and $q$.
\end{prop}

\begin{pf}
  Denote by $\k$ the adjoint trace field $\tr Ad \Gamma$, where
  $\Gamma=\CHL(G_{28},p,q)$. From the construction of the matrices
  $\widetilde{R}_j$, we get traces in $\Q(z,w)=\Q(\zeta_l)$, where
  $\zeta_d=e^{2\pi i/l}$. This implies $\k\subset
  \Q(\cos\frac{2\pi}{l})$.

  Since $\tr(R_1)=3+z$ and $\tr(R_3)=3+w$, we get
  $\cos(\frac{2\pi}{p}),\cos(\frac{2\pi}{q})\in\k$, which implies
  $\cos(\frac{2\pi}{l})\in\k$.
\end{pf}

\begin{prop}
  The CHL lattices derived from the group $G_{28}$ are all arithmetic.
\end{prop}

\begin{pf}
  We apply Theorem~\ref{thm:criterion} to the group generated by the
  matrices $\widetilde{R}_j$ given in
  equation~\eqref{eq:secondtry}. Note that the hypotheses of that
  theorem are satisfied, since the entries of $R_j$ are algebraic
  integers in the CM field $\Q(\zeta_l)$, and the adjoint trace field
  is equal to the maximal totally real subfield
  $\Q(\cos\frac{2\pi}{l})$ (as in Proposition~\ref{prop:pre}, $l$
  denotes the lcm of $p$ and $q$).

  For the list of pairs $(p,q)$ and the corresponding adjoint trace
  fields, see Table~\ref{tab:list}.
  The only cases that require work are those where the trace field is
  not $\Q$.  In each case, we need to compute the signature of
  non-trivial Galois conjugates of the Hermitian matrix of
  equation~\eqref{eq:form-g28}, and check that they are all definite
  (see Theorem~\ref{thm:criterion}).

  For $(p,q)=(2,5)$, up to complex conjugation, there is only one
  non-trivial Galois automorphism, given by $\zeta_{10}\mapsto
  \zeta_{10}^3$, which changes $\sqrt{5}$ to $-\sqrt{5}$. The
  Hermitian form~\eqref{eq:form-g28} has signature $(3,1)$ for $z=-1$,
  $w=\zeta_{10}^2=\zeta_5$, but it is positive definite for $z=-1$,
  $w=\zeta_{10}^6=\zeta_5^3$.

  For $(p,q)=(2,8)$, we need to consider the automorphism defined by
  $\zeta_{8}\mapsto \zeta_{8}^3$, which changes $\sqrt{2}$ to
  $-\sqrt{2}$.

  For $(p,q)=(2,12)$, $(3,4)$ or $(3,12)$, we need to consider
  $\zeta_{12}\mapsto \zeta_{12}^5$, which changes $\sqrt{3}$ to
  $-\sqrt{3}$.
\end{pf}

\subsection{Lattices derived from $G_{29}$} \label{sec:g29}

It follows from the results by Brou\'e, Malle, Rouquier~\cite{brmaro} and
Bessis and Michel~\cite{bessisannals},~\cite{bessismichel} that the
corresponding braid group is given by 
{\scriptsize
\begin{equation}\label{eq:braidgroup-g29}
\langle\,
r_1,r_2,r_3,r_4\, |\, 
\br_2(r_1,r_3), \br_2(r_1,r_4),
\br_3(r_1,r_2), \br_3(r_2,r_3), \br_3(r_3,r_4), 
\br_4(r_2,r_4), \br_4(r_3,r_2r_4)\,
\rangle.
\end{equation}
}

Couwenberg, Heckman and Looijenga show that there are representations
into $PU(3,1)$ with lattice image, that map every $r_j$ to a complex
reflection $R_j$ of angle $2\pi/p$, where $p$ is either $3$ or $4$. We
denote the corresponding groups by $\CHL(G_{29},p)$.

As before, we denote by $v_j$ a polar vector for the mirror of $R_j$
(this simply means that the mirror is the orthogonal complement of
$v_j$ with respect the Hermitian inner product). Note that these four
vectors must be linearly independent, because the group generated by
the $R_j$ must act irreducibly on $\C^4$. 

We take the vectors $v_j$ as the basis for $\C^4$, and because of
the braid relations $\br_3(R_j,R_{j+1})$, we can normalize them so
that the Hermitian form has the shape
\begin{equation}\label{eq:form_g29}
\left(\begin{matrix}
       1     &   \alpha   &      0     &   0\\
  \bar\alpha &     1      &   \alpha   & \beta\\
       0     & \bar\alpha &      1     & \alpha\\
       0     & \bar\beta  & \bar\alpha &   1
\end{matrix}\right)
\end{equation}
where $z=e^{2\pi i/p}$, $\alpha=1/(z-1)$ and $\beta$ is a complex
number to be determined.

By computing the matrices for the reflections $R_j$ and comparing the
(2,2)-entries of $(R_2R_4)^2$ and $(R_4R_2)^2$, it is easy to see that
the braid relation $\br_4(R_2,R_4)$ implies $|\beta|^2=0$ or
\begin{equation}\label{eq:g29-nsa}
|\beta|^2=\frac{2}{|z-1|^2}.
\end{equation}
The case $\beta=0$ is ruled out exactly as in
Proposition~\ref{prop:donotcommute}.

By computing the (2,2)-entry of $(R_3(R_2R_4))^2$ and
$((R_2R_4)R_3)^2$, we get the equation
$$
|\beta|^2z(z-1)-\bar\beta+\beta z^3=0,
$$
which, together with equation~\eqref{eq:g29-nsa}, implies
$$
\re(\beta z(z-1))=1,
$$
hence
$$
\beta=\frac{\mu}{z(z-1)},
$$
where $\mu=1\pm i$.

The corresponding matrices $R_j$ are given by 
{\tiny
\begin{eqnarray}\label{eq:matrices-g29}
&R_1=
  \left(\begin{matrix}
    z & 1 & 0 & 0\\
    0 & 1 & 0 & 0\\
    0 & 0 & 1 & 0\\
    0 & 0 & 0 & 1
  \end{matrix}\right),\quad 
R_2=
  \left(\begin{matrix}
    1  & 0 & 0 & 0\\
    -z & z & 1 & \mu\bar z\\
    0  & 0 & 1 & 0\\
    0  & 0 & 0 & 1
  \end{matrix}\right)
&R_3=
  \left(\begin{matrix}
    1 & 0  & 0  & 0\\
    0 & 1  & 0  & 0\\
    0 & -z & z  & 1\\
    0 &  0 & 0  & 1
  \end{matrix}\right),\quad
R_4=
  \left(\begin{matrix}
    1 & 0            & 0  & 0\\
    0 & 1            & 0  & 0\\
    0 & 0            & 1  & 0\\
    0 & -\bar\mu z^2 & -z & z
  \end{matrix}\right),
\end{eqnarray}
}
which preserve the Hermitian form
\begin{equation}\label{eq:g29-hf}
H = \left(\begin{matrix}
    1                  & \frac{1}{z-1}              & 0                  & 0\\
    \frac{1}{\bar z-1} & 1                          & \frac{1}{z-1}      & \frac{\mu\bar z}{z-1}\\
    0                  & \frac{1}{\bar z-1}         & 1                  & \frac{1}{z-1}\\
    0                  & \frac{\bar\mu z}{\bar z-1} & \frac{1}{\bar z-1} &  1
  \end{matrix}\right).
\end{equation}

One verifies that the Hermitian forms~\eqref{eq:g29-hf} corresponding
to the two choices $\mu=1\pm i$ both have signature $(3,1)$, for both
values $p=3$ and $p=4$. We now identify which choice corresponds to
the lattice $\CHL(G_{29},p)$.

The first remark is that in the case $p=3$, the two
matrices~\eqref{eq:g29-hf} are Galois conjugate. Indeed the
automorphism $\Q(\zeta_{12})$ that maps $\zeta_{12}$ to $\zeta_{12}^7$
fixes $\zeta_3$ while changing $\zeta_4$ to $\bar \zeta_4$.

\begin{prop}
  The group $\CHL(G_{29},p)$ for $p=3,4$ corresponds to choosing
  $\mu=1+i$ in formulas~\eqref{eq:matrices-g29} and~\eqref{eq:g29-hf}.
\end{prop}

\begin{pf}
  {\bf Let us first assume $p=3$.} We denote by $H^+$ (resp. $H^-$)
  the matrix $H$ of equation~\eqref{eq:g29-hf} for $\mu=1+i$
  (resp. $\mu=1-i$), and $z=e^{2\pi i/3}$. We also write $R_j^+$
  (resp. $R_j^-$) for the reflections preserving $H^+$ (resp. $H^-$),
  and finally we write $\Gamma^+$ (resp. $\Gamma^-$) for the group
  generated by the $R_j^+$ (resp. $R_j^-$).

  The lower right $3\times 3$ submatrix of $H^-$ gives a positive
  definite Hermitian form, so the subgroup $\Gamma^-_{234}$ of
  $\Gamma^-$ generated by $R_2^-,R_3^-,R_4^-$ has a fixed point inside
  the ball.

  One easily checks that $R_2^-R_3^-R_4^-$ is elliptic (the 1-eigenvector is
  negative for $H_-$), but has infinite order. The easiest way to
  check this is to consider the Galois conjugate group, where the
  corresponding matrix $R_2^+R_3^+R_4^+$ is loxodromic.

  This implies that the group corresponding to $\mu=1-i$ is not discrete.

  {\bf Let us now assume $p=4$.}  The argument uses a bit of
  CHL-theory, see~\cite{chl} or~\cite{chl3d_2}.  We consider the line
  $L=m_2\cap m_3\cap m_4$ which is the intersection of three mirrors
  of the reflections $S_2$, $S_3$ and $S_4$ in $G_{29}$ (as
  previously, $S_j$ and $R_j$ correspond to the same braid
  element). Its fixed point stabilizer in $G_{29}$ is isomorphic to
  the imprimitive Shephard-Todd group $G(4,4,3)$, which has order 96,
  and has 12 mirrors of reflections. This implies that $L$ is
  contained in 12 mirrors of the arrangement for $G_{29}$.

  For the group $\CHL(G_{29},4)$, the parameter $\kappa_L$ is given by
  $\kappa_L=\frac{12}{{\rm codim}
    L}(1-\frac{2}{p})=4(1-\frac{2}{4})=2$ (see p.~88
  of~\cite{chl}). The fact that $\kappa_L>1$ implies that the subgroup
  of $\CHL(G_{29},4)$ generated by $R_2,R_3$ and $R_4$ must preserve a
  totally geodesic copy of $H^2_\C$, hence the restriction of the
  Hermitian form to the span of the last three standard basis vectors
  $e_2,e_3,e_4$ must have signature $(2,1)$.

  One easily checks that, among the matrices $H^+$ and $H^-$ in
  equation~\eqref{eq:g29-hf} by taking ($z=i$ and) $\mu=1\pm i$, only
  $H_+$ has a lower-right $3\times 3$ block of signature $(2,1)$.
\end{pf}

\begin{prop}
  The adjoint trace field of $\CHL(G_{29},p)$ is $\Q(\sqrt{3})$ for
  $p=3$, and $\Q$ for $p=4$.
\end{prop}

\begin{pf}
  Given the above matrices and the fact that $\mu=1\pm i$, we clearly
  have matrices with entries in $\Q(i,z)$.

  For $p=4$, $\Q(i,z)=\Q(i)$, so the adjoint trace field is $\Q$.  For
  $p=3$, $\Q(i,z)=\Q(\zeta_{12})$, since $z$ is a primitive cube root
  of unity. This shows that the adjoint trace field is contained in
  $\Q(\sqrt{3})$.

  To show the other inclusion, we compute $\tr(R_4R_3R_2)=1+iz$, which
  for $z=\frac{-1+i\sqrt{3}}{2}$ gives
  $|\tr(R_4R_3R_2)|^2=2+\sqrt{3}$.
\end{pf}

\begin{prop}\label{prop:g29NA}
  The lattice $\CHL(G_{29},p)$ is non-arithmetic for $p=3$, and
  arithmetic for $p=4$.
\end{prop}

\begin{pf}
  Note that all the entries of the matrices in
  equation~\eqref{eq:matrices-g29} are algebraic integers, and their
  entries are in a CM field with maximal totally real subfield equal
  to the adjoint trace field, so we can apply
  Theorem~\ref{thm:criterion}.

  The arithmeticity in the case $p=4$ is obvious, since the trace
  field is $\Q$ and there is no nontrivial Galois conjugate to
  consider.

  The non-arithmeticity in the case $p=3$ follows immediately from the
  above discussion. Indeed, one checks that the matrices obtained in
  equation~\eqref{eq:g29-hf} for $\mu=1+i$ and $\mu=1-i$ both have
  signature (3,1).  
\end{pf}

\subsection{Lattices derived from the group $G_{30}$} \label{sec:g30}

The computations are similar to those in
section~\ref{sec:g28}. Because of the braid relations
$$ 
\br_2(R_1,R_3), \br_2(R_1,R_4), \br_2(R_2,R_4), \br_3(R_1,R_2),\br_3(R_2,R_3), \br_5(R_3,R_4)
$$
we may take the invariant Hermitian form to be
\begin{equation}\label{eq:herm-g30}
\left(\begin{matrix}
  1 & \alpha & 0 & 0\\
  {\bar \alpha} & 1 & \alpha & 0\\
  0 & {\bar \alpha} & 1 & \beta\\
  0  &  0  &  {\bar \beta} & 1
\end{matrix}\right),
\end{equation}
where
$$
\alpha = \frac{1}{z-1},\quad  \beta = \varphi\frac{1}{z-1},
$$ 
where $\varphi=\frac{1\pm\sqrt{5}}{2}$ (one could use
$\frac{-1\pm\sqrt{5}}{2}$, but this would give a conjugate subgroup of
$GL(4,\C)$).

In order to get a form with signature $(3,1)$ (at least when $p=3$ or
$5$), we need to take $\varphi=\frac{1+\sqrt{5}}{2}$, which we will
do in the sequel.

\begin{prop}\label{prop:adj30}
  The adjoint trace field of $\CHL(G_{30},p)$ is equal to
  $\Q(\cos{2\pi}{p},\varphi)$.
\end{prop}

\begin{pf}
  We denote by $\k$ the adjoint trace field, and by $\l$ the field
  $\Q(\cos\frac{2\pi}{p},\varphi)$.  Since the Hermitian
  form~\eqref{eq:herm-g30} is defined over $\Q(z,\varphi)$, we have
  $\k\subset\l$.
  
  On the other hand, we have $\tr(R_1)=3+z$ and
  $\tr(R_1R_2R_3R_4)=z(1-\varphi)$. This implies
  $|\tr(R_1)|^2=13+6\cos\frac{2\pi}{p}$ and
  $|\tr(R_1R_2R_3R_4)|^2=2-\varphi$, so $\l\subset\k$.
\end{pf}

From the Hermitian form, one can compute the matrices for $R_j$, which
has mirror polar to $e_j$ (the $j$-th vector in the standard basis for
$\C^4$), and multiplier $z=e^{2\pi i/p}$. We get
{\tiny
\begin{eqnarray}\label{eq:matrices-g30-1}
    R_1=\left(\begin{matrix}
      z & 1 & 0 & 0\\
      0 & 1 & 0 & 0\\
      0 & 0 & 1 & 0\\
      0 & 0 & 0 & 1
    \end{matrix}\right),\ 
    R_2=\left(\begin{matrix}
      1  & 0  & 0 & 0\\
      -z & z  & 1 & 0\\
      0  & 0  & 1 & 0\\
      0  & 0  & 0 & 1
    \end{matrix}\right),
    R_3=\left(\begin{matrix}
      1 & 0  & 0 & 0\\
      0 & 1  & 0 & 0\\
      0 & -z & z & \varphi\\
      0 & 0  & 0 & 1
    \end{matrix}\right),\ 
    R_4=\left(\begin{matrix}
      1 & 0 & 0          & 0\\
      0 & 1 & 0          & 0\\
      0 & 0 & 1          & 0\\
      0 & 0 & -\varphi z & \omega
    \end{matrix}\right).
\end{eqnarray}
}
Note that these matrices have entries in $\Z[z,\varphi]$.

\begin{prop} \label{prop:g30}
  The group $\mathcal{C}(G_{30},p)$ is a lattice if $p=3,5$, and in
  both cases it is cocompact. Both groups are arithmetic, with $\Q(\tr
  Ad \Gamma)=\Q(\sqrt{5})$.
\end{prop}

\begin{pf}
  The entries of the matrices in equation~\eqref{eq:matrices-g30-1}
  are algebraic integers in the CM field $\Q(z,\varphi)$, which is
  $\Q(i\sqrt{3},\sqrt{5})$ for $p=3$ and $\Q(\zeta_{5})$ for $p=5$.
  Both these fields have maximal totally real subfield equal to
  $\Q(\sqrt{5})$, which is equal to the adjoint trace field.

  In order to show arithmeticity, by Theorem~\ref{thm:criterion}, we
  need to show that non-trivial Galois conjugates of the Hermitian
  form~\eqref{eq:herm-g30} are definite.

  For $p=3$, up to complex conjugation, the only non-trivial Galois
  automorphism is given by $\sqrt{5}\mapsto -\sqrt{5}$ (and we may
  assume $i\sqrt{3}$ is left unchanged). The signature of the
  Hermitian form~\eqref{eq:herm-g30} is $(4,0)$ for
  $z=\frac{-1+i\sqrt{3}}{2}$ and $\varphi=\frac{1-\sqrt{5}}{2}$.

  For $p=5$, there is only one automorphism to consider (up to complex
  conjugation), given by $\zeta_5\mapsto\zeta_5^2$, which again
  changes $\frac{1+\sqrt{5}}{2}$ into $\frac{1-\sqrt{5}}{2}$.
  Taking $z=e^{4\pi i/5}$ and $\varphi = \frac{1-\sqrt{5}}{2}$ in the
  matrix~\eqref{eq:herm-g30} gives signature $(4,0)$.
\end{pf}

\subsection{Lattices derived from the group $G_{31}$} \label{sec:g31}

There are two groups in the CHL list, corresponding to $p=3$ and
$p=5$. These are a bit more difficult computationally, but not
conceptually.

The initial difficulty is that the corresponding Shephard-Todd is not
well-generated, i.e. it requires five generators (and not four as one
may expect from the dimension). We first parametrize quadruples of
reflections that satisfy the braid relations not involving $R_5$, and
express them in terms of off-diagonal entries of the Hermitian form,
see $\alpha$, $\beta$ in~\eqref{eq:H31}. Next, we use restrictions on
$\alpha$ and $\beta$ that come from the the existence of a 5-th
reflection that satisfies the appropriate relations with the first 4
reflections (this is expressed in terms of the parameters $\alpha$,
$\beta$ and the coordinates of a suitably normalized polar vector for
the reflection $R_5$, see the parameter $z$ below).

According to~\cite{bessismichel}, the group is generated by
reflections $R_1,\dots R_5$ that satisfy
\begin{equation}\label{eq:rels_g31}
\begin{array}{c}
\br_3(R_1,R_2), \br_3(R_2,R_5), \br_3(R_5,R_3), \br_3(R_3,R_4),\\
\br_2(R_2,R_4), \br_2(R_1,R_3), \br_2(R_2,R_3),\\
R_5R_4R_1=R_4R_1R_5=R_1R_5R_4.
\end{array}
\end{equation}
We denote by $v_j$ a polar vector to the mirror of $R_j$.  Note that
the last relation implies that the polar vectors $v_1,v_4,v_5$ are
linearly dependent. 

Since the action of the group generated by all the $R_j$ must be
irreducible on $\C^4$, the vectors $v_1,v_2,v_3,v_4$ must be linearly
independent. We write the Hermitian form in the corresponding
basis. The right angles coming from the above commutation relations
imply that we may assume the corresponding Hermitian matrix has the
form
\begin{equation}
H = \left(\begin{matrix}
  1 & \alpha & 0 & \beta\\
  \alpha & 1 & 0 & 0\\
  0 & 0 & 1 & \alpha\\
  \beta & 0 & \alpha & 1
\end{matrix}\right),\label{eq:H31}
\end{equation}
and we can choose $\alpha=1/(z-1)$ as in the previous sections.

We write $v_5$ for a polar vector to the mirror of $R_5$. Because of
the linear dependence between $v_1,v_4$ and $v_5$, we can write
$v_5=(x_1,0,0,x_4)$. If $R_1$ and $R_5$ have the same mirror, then
they coincide, and this would imply that $R_1,R_4$ commute, in which
case the action cannot be irreducible on $\C^4$.

Hence we must have $x_1\neq 0$, and we can take $w=(1,0,0,\mu)$ for
some $\mu\in\C$; by a similar reasoning, we must have $\mu\neq 0$.
Writing out matrices for $R_1$, $R_4$ and $R_5$ in terms of the
parameters $z,\beta,\mu$, and comparing the $(1,3)$-entries of
$R_5R_4R_1$ and $R_4R_1R_5$, we must have
$$
  |\mu|=1.
$$ 
Adjusting $v_4$ (and $v_3$ as well, since we want to keep the $\langle
v_3,v_4\rangle$ unchanged) by multiplying them by a complex number of
modulus one, we may assume further that $\mu=1$, i.e. we assume
$$
  v_5=(1,0,0,1).
$$

Comparing the $(1,2)$-entries of $R_5R_4R_1$ and $R_4R_1R_5$, we
now get the equation
\begin{equation}\label{eq:rel-541}
   |\beta|^2(1-z)+\beta-\bar\beta z=0.
\end{equation}

It is fairly easy to see that the braid relation $\br_3(R_2,R_5)$ is
equivalent to the relation
\begin{equation}\label{eq:rel-br3}
  \beta+\bar\beta+1=0.
\end{equation}

Combining equations~\eqref{eq:rel-541} and~\eqref{eq:rel-br3}, we get
\begin{equation}
  \beta=\frac{1}{\bar r-1},
\end{equation}
where $r=\sqrt{z}$ denotes one of the two complex numbers whose square
is $z$.

One then verifies (most reasonably with a computational software) that
all the relations in equation~\eqref{eq:rels_g31} are satisfied when we
take any of these two values of $\beta$. 

\begin{prop}
  The lattices $\CHL(G_{31},p)$ for $p=3,5$ are obtained by taking
  $z=e^{2\pi i/p}$, $r=e^{\pi i/p}$ in the Hermitian form $H$.
\end{prop}

\begin{pf}
  The fact that we take $z=e^{2\pi i/p}$ is from CHL theory. The two
  square-roots of $z$ are $\pm e^{\pi i/p}$. For $p=3$, $r=-e^{\pi
    i/3}$ gives a form of signature $(4,0)$. 
  
  The case $p=5$ is a bit more difficult, since both values $r=\pm
  e^{\pi i/5}$ give a form of signature $(3,1)$.

  In order to rule out $r=-e^{\pi i/5}$, we consider the $3\times 3$
  submatrix obtained from $H$ by removing the third row and column
  from $H$; for $r=-e^{\pi i/5}$, this submatrix gives a degenerate
  Hermitian form, which implies that $R_1, R_2, R_4$ have a common
  global fixed point at infinity (explicitly, this fixed point can be
  obtained by computing $e_1^\perp\cap e_2^\perp\cap e_4^\perp$ using
  linear algebra). One easily verifies that $R_1R_2R_4$ is a parabolic
  element (its eigenspace for the eigenvalue 1 is only 1-dimensional,
  but it is a double root of its characteristic polynomial).

  This rules out $r=-e^{\pi i/5}$, since the lattice $\CHL(G_{31},5)$
  is cocompact (see p.~160 of~\cite{chl}, where this group appears in
  boldface).
\end{pf}

We get the following matrices
{\small
\begin{eqnarray}\label{eq:matrices-g31-1}
  R_1=\left(\begin{matrix}
    z & 1 & 0 & -r-z\\
    0 & 1 & 0 &  0\\
    0 & 0 & 1 &  0\\
    0 & 0 & 0 &  1
  \end{matrix}\right),\quad
  R_2=\left(\begin{matrix}
    1  & 0 & 0 & 0\\
    -z & z & 0 & 0\\
    0  & 0 & 1 & 0\\
    0  & 0 & 0 & 1
  \end{matrix}\right),\quad
  R_3=\left(\begin{matrix}
    1 & 0 & 0 & 0\\
    0 & 1 & 0 & 0\\
    0 & 0 & z & 1\\
    0 & 0 & 0 & 1
  \end{matrix}\right)\\
  R_4=\left(\begin{matrix}\label{eq:matrices-g31-2}
    1      & 0 & 0  & 0 \\
    0      & 1 & 0  & 0\\
    0      & 0 & 1  & 0\\
   r+1     & 0 & -z & z
  \end{matrix}\right),\quad
  R_5=\left(\begin{matrix}
    z+r+1 & 1 & -z & -r-1\\
    0     & 1 & 0  & 0\\
    0     & 0 & 1  & 0\\
    z+r   & 1 & -z & -r
  \end{matrix}\right),
\end{eqnarray}
}
where $z=e^{2\pi i/p}$ and $r=e^{\pi i/p}$.

We only consider the cases where $p=3$ or $5$, i.e. $z$ is a root of
unity of odd order, so the corresponding cyclotomic field $\Q(z)$
contains the square-roots of $z=e^{2\pi i/p}$ (if $p=2m-1$, the
square-roots of $z$ are $\pm z^m$).  In particular we get that, in
both cases, the adjoint trace field is contained in
$\Q(\cos\frac{2\pi}{p})$ and, as before, the fact that $\tr(R_1)=3+z$
implies:
\begin{prop}
  For $p=3$ or $5$, the adjoint trace field of $\CHL(G_{31},p)$ is
  $\Q(\cos\frac{2\pi}{p})$.
\end{prop}

\begin{prop}
  For $p=3$ and $5$, the lattice $\CHL(G_{31},p)$ is arithmetic.
\end{prop}

\begin{pf}
  The entries of the matrices in
  equations~\eqref{eq:matrices-g31-1},~\eqref{eq:matrices-g31-2} are
  algebraic integers in the cyclotomic field $\Q(z)$ where $z=e^{2\pi
    i/p}$, and we can apply the arithmeticity criterion of
  Theorem~\ref{thm:criterion}.
  
  For $p=3$, there is nothing to check, since the adjoint trace field
  is $\Q$ and there is no nontrivial Galois conjugate.

  For $p=5$, the above Hermitian form reads
  \begin{equation}
    H = \left(\begin{matrix}
      1                       & \frac{1}{z-1}     & 0                     & -\frac{1}{z^2+1}\\
      \frac{1}{\bar z-1}      & 1                 & 0                     & 0\\
      0                       & 0                 & 1                     & \frac{1}{z-1}\\
      -\frac{1}{\bar z^2+1}   & 0                 & \frac{1}{\bar z-1}    & 1
     \end{matrix}\right),
  \end{equation}
  where $z=e^{2\pi i/5}$, and this Hermitian form has signature $(3,1)$.

  When changing $z$ to $z^2$ (which changes $\sqrt{5}$ to
  $-\sqrt{5}$), the form $H$ becomes definite, so
  Theorem~\ref{thm:criterion} says that $\CHL(G_{31},5)$ is
  arithmetic.
\end{pf}

\subsection{Lattices derived from the group $G_{33}$} \label{sec:g33}

Note that the group $G_{33}$ acts on $\C^5$ (hence the corresponding
CHL lattice acts on $H^4_{\mathbb{C}}$). More generally, here and in
the following few sections, the dimension of the relevant complex
hyperbolic space $H^n_{\mathbb{C}}$ can be read off the size
$(n+1)\times(n+1)$ of the Hermitian matrix.

The group $\CHL(G_{33},3)$ is generated by 5 reflections
$R_1,\dots,R_5$ of order 3. As before, we write the Hermitian form $H$
in the basis of $\C^5$ given by vectors polar to the mirrors of the
reflections $R_j$. By suitably rescaling these vectors, we may assume
$H$ has the shape
\begin{equation}\label{eq:H33}
H=\left(\begin{matrix}
       1     &        \alpha         &      0     &       0       &   0     \\
  \bar\alpha &           1           &   \alpha   & \lambda\alpha &   0     \\
       0     &      \bar\alpha       &      1     &    \alpha     &   0     \\
       0     & \bar\lambda\bar\alpha & \bar\alpha &       1       & \alpha  \\
       0     &           0           &      0     & \bar\alpha    &   1
\end{matrix}\right),
\end{equation}
where $\alpha=1/(z-1)$ as above, and $\lambda\in \C$ is to be
determined later.

The braid relation $\br_3(R_2,R_4)$ is equivalent to
$|\lambda|=1$. According to Bessis and Michel, we must have
\begin{equation}\label{eq:g33-relation}
  (R_2R_4R_3)^2=(R_4R_3R_2)^2=(R_3R_2R_4)^2.
\end{equation}

One can easily write down matrices for the $R_j$ in terms of $\lambda$
by using equation~\eqref{eq:refl} for $\langle v,w\rangle = w^*Hv$. By
comparing the (3,1) entries of $(R_2R_4R_3)^2$ and $(R_4R_3R_2)^2$, we
see that the relation~\eqref{eq:g33-relation} implies $\Re(\lambda
z)=\frac{1}{2}$, which gives $\lambda=\frac{1\pm \sqrt{3}}{2}\bar z$.

It is then easy to verify (most conveniently with some computer
algebra system) that both values of $\lambda$ make
relation~\eqref{eq:g33-relation} hold. 

We will only consider the case $p=3$. We write $\omega=e^{2\pi
  i/3}=\frac{-1+i\sqrt{3}}{2}$, and take $z=\omega$ in the above
matrix. The two values of $\lambda$ are then given by $-\omega$ and
$-\bar\omega$.

One verifies that, in that case, the Hermitian form in
equation~\eqref{eq:H33} is degenerate for $\lambda=-\omega$, whereas
is has signature $(4,1)$ for $\lambda=-\bar\omega$. Hence we have the
following.
\begin{prop}\label{prop:G33}
  The group $\CHL(G_{33},3)$ is isomorphic to the group generated by
  $R_j$, $j=1,\dots,5$ as above, preserving the Hermitian
  form~\eqref{eq:H33} for $z=\omega$ and $\lambda=-\bar\omega$.
\end{prop}
We write $w_j$ for the $j$-th row of $R_j$. Note that $R_j-Id$ has
only one non-zero row, so in order to describe $R_j$, it is enough to
list $w_j$, which is what we do in equation~\eqref{eq:rowsRj-33}.
\begin{equation}\label{eq:rowsRj-33}
  \begin{array}{l}
    w_1=(\omega,1,0,0,0)\\
    w_2=(-\omega,\omega,1,-\omega,0)\\
    w_3=(0,-\omega,\omega,1,0)\\
    w_4=(0,1,-\omega,\omega,1)\\
    w_5=(0,0,0,-\omega,\omega)
  \end{array}
  \quad
  H=\left(\begin{matrix}
  1                       & \frac{1}{\omega-1}   & 0                      & 0                      & 0 \\
  \frac{1}{\bar \omega-1} & 1                    & \frac{1}{\omega-1}     & \frac{1}{\bar\omega-1} & 0\\
  0                       & \frac{1}{\bar\omega-1}   & 1                      & \frac{1}{\omega-1}     & 0\\
  0                       & \frac{1}{\omega-1}   & \frac{1}{\bar\omega-1} & 1                      & \frac{1}{\omega-1}\\
  0                       & 0                    & 0                      & \frac{1}{\bar\omega-1} & 1
\end{matrix}\right).
\end{equation}

Equation~\eqref{eq:rowsRj-33} makes it clear that the matrices $R_j$
have algebraic integer entries in $Q(\omega)$, which implies $\Q(\tr
Ad\Gamma)=\Q$, so we have the following.
\begin{prop}
  The lattice $\CHL(G_{33},3)$ is arithmetic with adjoint trace field $\Q$.
\end{prop}

\subsection{The group $G_{34}$} \label{sec:g34}

The group $\CHL(G_{34},3)$ is generated by 6 reflections
$R_1,\dots,R_6$, the braid group is the same as the previous one, with
one extra generator that commutes with the first four, and braids with
length 3 with the fifth.

The same computations as in section~\ref{sec:g33} show that we can take $H$ to be
\begin{equation}\label{eq:H34}
H=\left(\begin{matrix}
           1             &          \frac{1}{\omega-1}              &         0              &          0               &        0               &       0 \\
  \frac{1}{\bar\omega-1} &                1                         &   \frac{1}{\omega-1}   & \frac{\lambda}{\omega-1} &        0               &       0\\
           0             &        \frac{1}{\bar\omega-1}            &         1              &    \frac{1}{\omega-1}    &        0               &       0\\
           0             & \frac{\overline{\lambda}}{\bar \omega-1} & \frac{1}{\bar\omega-1} &          1               & \frac{1}{\omega-1}     &       0\\
           0             &                0                         &         0              & \frac{1}{\omega-1}       &        1               & \frac{1}{\omega-1}\\
           0             &                0                         &         0              &          0               & \frac{1}{\bar\omega-1} &       1
\end{matrix}\right),
\end{equation}
where again $\omega=\frac{-1+i\sqrt{3}}{2}$ and $\lambda$ is either
$-\bar\omega$ or $-\omega$. One readily checks that both Hermitian
forms have signature $(5,1)$ so it is not clear which group
corresponds to $\CHL(G_{34},3)$.

We call $H^+$ (resp. $H^-$) the Hermitian form corresponding to
$\lambda=\frac{1+i\sqrt{3}}{2}$
(resp. $\lambda=\frac{1-i\sqrt{3}}{2}$). Note that the upper left
$4\times 4$ submatrices $K^\pm$ of $H^\pm$ do not have the same
signature, namely $K^+$ has signature $(3,1)$, whereas $K^-$ is
degenerate.

Using a bit of CHL theory (see the argument below), the last
observation implies the following.
\begin{prop}\label{prop:G34}
  The group $\CHL(G_{34},3)$ is conjugate to the group generated by
  the complex reflections $R_j$ with multiplier $\omega$ and polar
  vectors given by the standard basis vectors $e_j$ of $\C^6$, and
  Hermitian form $H^+$.
\end{prop}
\begin{pf}
  Denote by $S_j$, $j=1,\dots,6$ the reflections generating $G_{34}$
  (with the same numbering as in Figure~\ref{fig:cox}(h) on
  page~\pageref{fig:cox}); moreover, we denote by $R_j$ the
  corresponding reflections in $\CHL(G_{34},3)$ (i.e. $R_j$ and $S_j$
  are the images of the same element $r_j$ of the relevant braid
  group, in the Bessis-Michel presentations).

  %% the numbers below probably need to be modified, they were for  
%%G_{33]...  
  It is quite clear from the Coxeter diagrams in Figure~\ref{fig:cox}
  that $S_1,S_2,S_3,S_4$ and $S_5$ generate a group isomorphic to
  $G_{33}$ (this is a group of order 51840).
  The 1-dimensional intersection of their mirrors in $\C^6$ is
  contained in 45 mirrors of reflections in $G_{34}$ (this is the
  number of mirrors of reflections in $G_{33}$, see p.~302
  in~\cite{shephardtodd}).

  The weight attached to the stratum $L=L_{12345}$, i.e. the
  intersection of the mirrors of $S_1,\dots,S_5$ is
  $$
  \kappa_L=\frac{45}{5}(1-\frac{2}{3})=3>1,
  $$ 
  where the denominator 5 comes from the codimension of this stratum
  (see p.~96 of~\cite{chl}). Since $\kappa_L>1$,
  CHL theory predicts that the mirrors of
  $R_1,R_2,R_3,R_4$ and $R_5$ should be orthogonal to a common complex
  hyperbolic totally geodesic copy of $H^4_\mathbb{C}$ in
  $H^5_\mathbb{C}$. 

  Such a totally geodesic copy it given by the orthogonal complement
  $v^\perp$ of a vector $v$ with $\langle v,v\rangle>0$, so the
  restriction of the Hermitian form to the complex span of
  $e_1,\dots,e_5$ must have signature $(4,1)$, so we must use $H^+$
  (and not $H^-$).
\end{pf}

For completeness, as in section~\ref{sec:g33}, we describe the
(non-obvious rows of) the matrices $R_j$, $j=1,\dots,6$ in
equation~\eqref{eq:rowsRj-34} (recall that $w_j$ is the $j$-th row of
$R_j$).
\begin{equation}\label{eq:rowsRj-34}
  H=\left(\begin{matrix}
  1                       & \frac{1}{\omega-1}   & 0                      & 0                      & 0                      & 0\\
  \frac{1}{\bar \omega-1} & 1                    & \frac{1}{\omega-1}     & \frac{1}{\bar\omega-1} & 0                      & 0\\
  0                       & \frac{1}{\bar\omega-1}   & 1                      & \frac{1}{\omega-1}     & 0                      & 0\\
  0                       & \frac{1}{\omega-1}   & \frac{1}{\bar\omega-1} & 1                      & \frac{1}{\omega-1}     & 0\\
  0                       & 0                    & 0                      & \frac{1}{\bar\omega-1} &           1            &  \frac{1}{\omega-1}\\
  0                       & 0                    & 0                      & 0                      & \frac{1}{\bar\omega-1} &     1
\end{matrix}\right).
\quad 
 \begin{array}{l}
    w_1=(\omega,1,0,0,0,0)\\
    w_2=(-\omega,\omega,1,-\omega,0,0)\\
    w_3=(0,-\omega,\omega,1,0,0)\\
    w_4=(0,1,-\omega,\omega,1,0)\\
    w_5=(0,0,0,-\omega,\omega,1)\\
    w_6=(0,0,0,0,-\omega,\omega)
  \end{array}
\end{equation}

The matrices $R_1,\dots,R_6$ have entries in $\Z[\omega]$, which
implies once again that $\Q(\tr Ad\Gamma)=\Q$, so we get the
following.
\begin{prop}
  The lattice $\CHL(G_{34},3)$ is arithmetic with adjoint trace field $\Q$.
\end{prop}

\subsection{Lattices derived from the group $G_{35}$} \label{sec:g35}

In this case, we can write the Hermitian matrix as
$$
\left(\begin{matrix}
  1 & \alpha& 0 & 0 & 0 & 0\\
  \bar \alpha& 1 & \alpha& 0 & 0 & 0\\
  0 & \bar \alpha& 1 & \alpha& 0 & \alpha\\
  0 & 0 & \bar \alpha& 1 & \alpha& 0\\
  0 & 0 & 0 & \bar \alpha& 1 & 0\\
  0 & 0 & \bar \alpha& 0 & 0 & 1
\end{matrix}\right).
$$
where $\alpha=\frac{1}{z-1}$, $z=e^{2\pi i/p}$.
The corresponding reflections have entries in $\Z[z]$, hence we have the following.
\begin{prop}
  The lattices $\CHL(G_{35},3)$ and $\CHL(G_{35},4)$ are both arithmetic with adjoint trace field $\Q$.
\end{prop}

\subsection{The group $G_{36}$} \label{sec:g36}

Note that $G_{36}$ acts on $\C^7$, the corresponding CHL lattice acts
on $H^6_{\mathbb{C}}$. We can write the Hermitian matrix as
$$
\left(\begin{matrix}
  1 & \alpha & 0 & 0 & 0 & 0 & 0\\
  \bar \alpha & 1 & \alpha & 0 & 0 & 0 & 0\\
  0 & \bar \alpha & 1 & \alpha & 0 & \alpha & 0\\
  0 & 0 & \bar \alpha & 1 & \alpha & 0 & 0\\
  0 & 0 & 0 & \bar \alpha & 1 & 0 & \alpha\\
  0 & 0 & \bar \alpha & 0 & 0 & 1 & 0\\
  0 & 0 & 0 & 0 & \bar \alpha & 0 & 1
\end{matrix}\right).
$$ 
where $\alpha=\frac{1}{\omega-1}$. This gives matrices with entries in
$\Z[\omega]$, so we have:
\begin{prop}
  The lattices $\CHL(G_{36},3)$ is arithmetic with adjoint trace field $\Q$.
\end{prop}

\subsection{Lattices derived from the group $G_{37}$} \label{sec:g37}

We can write the Hermitian matrix as
$$
\left(\begin{matrix}
  1 & \alpha& 0 & 0 & 0 & 0 & 0 & 0\\
  \bar \alpha& 1 & \alpha& 0 & 0 & 0 & 0 & 0\\
  0 & \bar \alpha& 1 & \alpha& 0 & \alpha& 0 & 0\\
  0 & 0 & \bar \alpha& 1 & \alpha& 0 & 0 & 0\\
  0 & 0 & 0 & \bar \alpha& 1 & 0 & \alpha& 0\\
  0 & 0 & \bar \alpha& 0 & 0 & 1 & 0 & 0\\
  0 & 0 & 0 & 0 & \bar \alpha& 0 & 1 & \alpha\\
  0 & 0 & 0 & 0 & 0 & 0 & \bar \alpha& 1
\end{matrix}\right).
$$
where $\alpha=\frac{1}{\omega-1}$. This gives matrices with entries in
$\Z[\omega]$, so we have:
\begin{prop}
  The lattices $\CHL(G_{37},3)$ is arithmetic with adjoint trace field $\Q$.
\end{prop}

\section{The proof of Theorem~\ref{thm:new3d}}\label{sec:new3d}

The fact that $\CHL(G_{29},3)$ is not cocompact is mentioned in the
tables in~\cite{chl}. The adjoint trace field was determined in
section~\ref{sec:arithmeticity}, where we also proved
non-arithmeticity, see Proposition~\ref{prop:g29NA}. The only thing
that is left to prove is the fact that it is not commensurable to the
Deligne-Mostow group $\Gamma_\mu$ with $\mu=(3,3,3,3,5,7)/12$.

This is not obvious, since both groups have the same rough
commensurability invariants (both are non-uniform, have
non-arithmeticity index one, and they have the same adjoint trace
field).

We will argue by comparing the cusps in both groups. It is known that
both groups have a single orbit of cusps, but we will show that the
corresponding cusps are not commensurable.  In
section~\ref{sec:cuspg29} and~\ref{sec:cuspdm} we will describe their
respective cusps, and in~\ref{sec:incomm} we will use this to show
that they are incommensurable (see Proposition~\ref{prop:incomm}).
We start with a general descriptions of cusps in 3-dimensional complex
hyperbolic space (most of this can be found in chapter~4
of~\cite{goldman}).

\subsection{Cusps and the Heisenberg group} \label{sec:cusps_general}

When studying a cusp, we will write the Hermitian form in block form as
\begin{equation}\label{eq:herm_cusp}
  H = \left(\begin{array}{c|c|c}
    0 & 0 & 1\\
 \hline
    0 & K & 0\\
 \hline
    1 & 0 & 0
  \end{array}\right),
\end{equation}
where $K$ is a positive definite $2\times 2$ Hermitian form.

It is easy to see that any parabolic transformation in $PU(H)$ fixing
$(1,0,0,0)$ can be written as
\begin{equation}\label{eq:parabolic}
  P(B,w,t) = \left(\begin{array}{c|c|c}
    1 & -w^*KB & -\frac{1}{2}w^*Kw+it\\
 \hline
    0 & B & w\\
 \hline
    0 & 0 & 1
  \end{array}\right),
\end{equation}
where $B\in U(K)$, $w\in\C^2$ and $t\in\R$.

Moreover, one has
$$
  P(B,w,t)P(B',w',t')=P(BB',Bw'+w,t+t'+\im\left( w'^* B^* K w \right)),
$$
and
$$
  P(B,w,t)^{-1}=P(B^{-1},-B^{-1}w,-t).
$$

In the special case $B=Id$, we get unipotent elements
$U(w,t)=P(Id,w,t)$ that satisfy
\begin{equation}\label{eq:heisenberg}
  U(w,t)U(w',t')=U(w+w',t+t'+\im\left( w'^* K w \right)), \quad U(w,t)^{-1}=U(-w,-t).
\end{equation}
The corresponding group law $(w,t)\star(w',t)=(w+w',t+t'+\im\left(
w'^* K w \right))$ is called the (5-dimensional) Heisenberg group law,
and the corresponding isometry of complex hyperbolic space is called a
Heisenberg translation. When $w=0$, the corresponding translation is
called a \emph{vertical} translation, and these are precisely the
elements that are central in the Heisenberg group.

From equation~\eqref{eq:heisenberg}, one easily checks that the
commutator of $U(w,t)$ and $U(w',t')$ is given by
$$
  [U(w,t),U(w',t')]=U(0,2\im\left( w'^*Kw \right)),
$$
which is a vertical translation.

The parabolic stabilizer of $(1,0,0,0)$ has a projection onto the
complex unitary affine group $\C^2\ltimes U(K)$, given by keeping only
the lower-right $3\times 3$ block of the above matrices $P(B,w,t)$,
the kernel consisting of the vertical translations in the group.

We will use this description to describe a cusp $\Gamma_\infty$ of a
non-cocompact lattice $\Gamma$. We may assume that the corresponding
ideal fixed point is given by $(1,0,0,0)$, and that the Hermitian form
is as in equation~\eqref{eq:herm_cusp}. Since a discrete group cannot
have both parabolic and loxodromic elements fixing the same ideal
point, the cusps will only contain parabolic elements.

The projection of $\Gamma_\infty$ to the affine group $\C^2\ltimes
U(K)$ is then a complex crystallographic group, i.e. it must act
cocompactly (which amounts to requiring that its translation subgroup
has rank 4), and the vertical part is then an infinite cyclic group,
commensurable to the one obtained by taking a single non-trivial
commutator of Heisenberg translations.

Later in the paper, we will need to understand how the horizontal
($w$) and vertical ($t$) components of $U(w,t)$ behave under isometric
changes of coordinates. If $Q$ is a general isometry of the form $H$
in equation~\eqref{eq:herm_cusp}, then it can be written as
\begin{equation}\label{eq:isometry}
  Q = \left(\begin{array}{c|c|c}
    \alpha & -\alpha v^*KC & -\frac{1}{2}\alpha w^*Kw+is\\
 \hline
    0 & C & v\\
 \hline
    0 & 0 & \frac{1}{\alpha}
  \end{array}\right),
\end{equation}
where $C\in U(K)$, $\alpha,t\in\R$ and $v\in\C^2$. One verifies by
direct computation that 
\begin{equation}\label{eq:key}
  Q U(w,t) Q^{-1} = U(\alpha C w,\alpha^2(t+2\Im(w^*C^*Kv))).
\end{equation}
In particular, the square norm (with respect to $K$) of the horizontal
part gets multiplied by $\alpha^2$. In the special case $w=0$, we have
$Q U(0,t) Q^{-1} = U(0,\alpha^2 t)$ for every $t\in\R$, so the
vertical component of vertical translations also gets multiplied by
$\alpha^2$.

\subsection{The cusp of the Deligne-Mostow non-arithmetic lattice in $PU(3,1)$.} \label{sec:cuspdm}

We now review some facts about the cusp of the group
$\Gamma=\Gamma_{\mu,\Sigma}$ where $\mu=(3,3,3,3,5,7)/12$, and
$\Sigma\simeq S_4$ permutes the first four weights. It is well known
from Deligne-Mostow theory that the quotient has a single cusp, and
that the corresponding cusp stabilizer can be described using
hypergeometric functions for weights $(3,3,3,3,5+7)/12=(1,1,1,1,4)/4$,
and this monodromy group is described in section \S 15.20
of~\cite{delignemostowbook} (it corresponds to a parabolic case in
Deligne-Mostow).

We give another argument, based on the description of
$\Gamma=\Gamma_{\mu,\Sigma}$ as $\CHL(B_4,3,4)$ (see section~4
in~\cite{chl3d_2}). The point is that the relevant braid group is
generated by three reflections $R_2,R_3,R_4$ that correspond to
half-twists between points with equal weights (these have order
$2(1-2\mu_1)^{-1}=4$) and a complex reflection $R_1$ corresponding to
a full-twist between a point with weight $\mu_4$ and one with weight
$\mu_5$ (this has order $(1-\mu_4-\mu_5)^{-1}=3$). 

We then have $\br(R_2,R_3)=\br(R_3,R_4)=3$, $\br(R_1,R_2)=4$, and
other pairs of reflections among the $R_j$ commute. These braid
relations determine the group up to conjugation, using the fact that
$R_1$ (resp. $R_j$ for $j=2,3,4$) has multiplier $e^{2\pi i/3}=\omega$
(resp. $e^{2\pi i/4}=i$).

We describe the group in coordinates that make the structure of
the cusp group visible. Consider the Hermitian form
\begin{equation}\label{eq:HcuspDM}
H=\left(\begin{matrix}
  0 &   0  &   0  & 1\\
  0 &   2  & -1-i & 0\\
  0 & -1+i &   2  & 0\\
  1 &   0  &   0  & 0
\end{matrix}\right),
\end{equation}
and the reflections $R_1=R_{v_1,\omega}$, $R_j=R_{v_j,i}$ for
$j=2,3,4$, where we use the notation in equation~\eqref{eq:refl}, and
\begin{equation}\label{eq:vlist}
v_1=(2,0,0,\zeta^2+\zeta+1),\quad 
v_2=(2,i-1,0,0),\quad
v_3=(0,0,1,0),\quad
v_4=(0,i+1,2,0).
\end{equation}
Here we denote by $\zeta=e^{2\pi i/12}$. One easily verifies that
these satisfy the correct braid relations, and the corresponding
complex Coxeter diagram has no loop, so we get:
\begin{prop}
  The group generated by $R_1,\dots,R_4$ is conjugate to
  $\CHL(B_4,3,4)$, which is also $\Gamma_{\mu,\Sigma}$ for
  $\mu=(3,3,3,3,5,7)/12$ and $\Sigma=S_4$.
\end{prop}
Note that the non-trivial Galois conjugate of $\CHL(B_4,3,4)$ is also
generated by complex reflections satifsying the same relations, but
the multipliers of the complex reflections are different; indeed, it
corresponds to the automorphism $\varphi$ of $\Q(\zeta)$ induced by
$\varphi(\zeta)=\zeta^5$ (or its complex conjugate), and one then has
$\varphi(i)=\varphi(\zeta^3)=\zeta^{15}=\zeta^3=i$, and similarly
$\varphi(\omega)=\varphi(\zeta^4)=\zeta^{20}=\bar\omega$.

We now study the structure of the cusp group using the notation in
section~\ref{sec:cusps_general}. It follows from CHL theory that
$\CHL(B_4,3,4)$ has a single cusp, represented by the group generated
by $R_2,R_3$ and $R_4$.

Indeed, the cusps of the quotient are in one-to-one correspondence
with the $G$-orbits of irreducible mirror intersections $L$ with
$\kappa_L=1$, in the notation of~\cite{chl}. Recall that
$$
\kappa_L=\frac{\sum_{H\supset L}\kappa_H}{{\rm codim} L},
$$
where the sum ranges over all mirrors in the arrangement.

We denote by $r_1,\dots,r_4$ the generators of $G=W(B_4)$, by $L_j$
the mirror of $r_j$. A list of representatives for all (non-trivial,
pairwise distinct) irreducible mirror intersections is given in
Table~\ref{tab:irred_b4} (compare with the tables in~\cite{chl3d_2}). In
the table, $L_{jk}$ stands for $L_j\cap L_k$, and $L_{jkl}$ for
$L_j\cap L_k\cap L_l$.
\begin{table}
\begin{tabular}{|c||c|c|c|c|c|c|}
\hline
$L$        &         $L_1$      &       $L_2$       & $L_{12}$                        & $L_{23}$                       & $L_{123}$                       & $L_{234}$\\
\hline
Type       &        $(1,0)$     &      $(0,1)$      &  $(2,2)$                        &  $(0,3)$                       &  $(3,6)$                        &  $(0,6)$\\
\hline
$\kappa_L$ &  $1-\frac{2}{p_1}$ & $1-\frac{2}{p_2}$ & $2-\frac{2}{p_1}-\frac{2}{p_2}$ & $\frac{3}{2}(1-\frac{2}{p_2})$ & $3-\frac{2}{p_1}-\frac{4}{p_2}$ & $2(1-\frac{2}{p_2})$\\
\hline
\end{tabular}
\caption{Irreducible mirror intersections for $W(B_4)$. The type of an
  irreducible mirror intersection $L$ is the number of mirrors in each
  $G$-orbit of mirrors containing $L$, for instance an $L$ of type
  (3,6) is contained in 3 mirrors in the $G$-orbit of $L_1$ and 6
  mirrors in the $G$-orbit of $L_2$.}\label{tab:irred_b4}
\end{table}
For the group $\CHL(B_4,3,4)$, we take $(p_1,p_2)=(3,4)$, and only one
irreducible mirror intersection gives $\kappa_L=1$, namely
$L=L_{234}$.

In other words, the lattice $\Gamma=\CHL(B_4,3,4)$ has a single
$\Gamma$-orbit of cusps, which is generated by $R_2,R_3$ and
$R_4$. Note also that the corresponding cusp group is a CHL group of
type $\CHL(A_3,4)$, and it fits in their framework as a parabolic case
(see section 5.3 of~\cite{chl}).

We now work out a detailed description of the cusp group. The vector
$e_1=(1,0,0,0)$ is a null vector for the Hermitian form $H$, and
$\langle e_1,v_2\rangle=\langle e_1,v_3\rangle=\langle
e_1,v_4\rangle=0$ (see the definition of the $v_j$ in
equation~\eqref{eq:vlist}). This implies that $e_1$ gives the global
ideal fixed point of the group generated by $R_2,R_3,R_4$.

We have chosen the Hermitian form so that the corresponding matrices
for $R_2,R_3,R_4$ have Gaussian integers entries, i.e. entries in
$\Z[i]$ (the entries of $R_1$ are not algebraic integers, but this is
irrelevant). In fact, the corresponding reflections read { \small
\begin{eqnarray}\label{eq:matrices-DMNA}
  R_2=\left(\begin{matrix}
    1  & 2 & -1-i & -1+i\\
    0  & i &   1  & -i\\
    0  & 0 &   1  & 0\\
    0  & 0 &   0  & 1
  \end{matrix}\right),
  R_3=\left(\begin{matrix}
    1 & 0  & 0 & 0\\
    0 & 1  & 0 & 0\\
    0 & -i & i & 1\\
    0 & 0  & 0 & 1
  \end{matrix}\right),
  R_4=\left(\begin{matrix}
    1  & 0 & 0  & 0 \\
    0  & 1 & -1 & 0\\
    0  & 0 & i  & 0\\
    0  & 0 & 0  & 1
  \end{matrix}\right)
\end{eqnarray}.
} 
As mentioned in the general description of cusps given in
section~\ref{sec:cusps_general}, the projection onto the complex
affine group $U(K)\ltimes \C^2$ is given by the lower-right $3\times
3$ submatrices of $R_2,R_3$ and $R_4$.

The linear parts are given by
\begin{equation}\label{eq:B}
B_2=\left(\begin{matrix}
  i & 1\\
  0 & 1
\end{matrix}\right),\quad
B_3=\left(\begin{matrix}
  1  & 0\\
  -i & i
\end{matrix}\right),\quad
B_4=\left(\begin{matrix}
  1 & -1\\
  0 & i
\end{matrix}\right),
\end{equation}
and the translation parts are given by
\begin{equation}\label{eq:w}
w_2=\left(\begin{matrix}-i\\0 \end{matrix}\right),
w_3=\left(\begin{matrix}0\\1 \end{matrix}\right),
w_4=\left(\begin{matrix}0\\0 \end{matrix}\right).
\end{equation}

\begin{prop}
  The matrices $B_2,B_3,B_4$ generate a copy of the Shephard-Todd
  group $G_8$, which has order 96 (and center of order 4).
\end{prop}
\begin{pf}
  The $B_j$ are complex reflections of order $4$, and
  $\br(R_2,R_3)=3$. This implies that $B_2$ and $B_3$ generate a copy
  of $G_8$, and one easily verifies that $B_4$ is in the group
  generated by $B_2$ and $B_3$ (in fact $B_4=B_3B_2B_3B_2B_3$).
\end{pf}
This gives the first part of the following proposition.
\begin{prop}
  The projection to the complex affine group $U(K)\ltimes \C^2$ has
  linear part $G_8$, and translation subgroup given by
  $\Z[i]\times\Z[i]$.
\end{prop}

\begin{pf}
It is clear that the translation subgroup is contained in
$\Z[i]\times\Z[i]$. We claim that it is precisely equal to it, which
can be checked by computing the following,
\begin{equation}\label{eq:tsllist}
  \begin{array}{l}
  T_{(1,0)} = A_2 A_3^2 A_4^{-1} A_3^2\\
  T_{(i,0)} = A_3^2 A_4^{-1} A_3^2 A_2\\
  T_{(0,1)} = A_2^{-1} A_3 A_2 A_4^{-1} A_3 A_2^{-1}\\
  T_{(0,i)} = A_2 A_3^{-1} A_4^{-1} A_3 A_2^{-1} A_3.
  \end{array}
\end{equation}
Here we denote by $T_w$ the vertical translation $U(Id,0,w)$.
\end{pf}

In particular, given the shape of the $2\times 2$ Hermitian matrix $K$
in equation~\eqref{eq:HcuspDM}, we get that $w^*Kw\in \Z$ for every
translation.

The following follows from the fact that the matrices in
equation~\eqref{eq:matrices-DMNA} have entries in $\Z[i]$.
\begin{prop}
  The vertical translations in the cusp have translation length in
  $\Z$.
\end{prop}
In other works, the vertical translations are of the form $P(Id,0,t)$
for some $t\in\Z$.

\subsection{The cusp of $\CHL(G_{29},3)$} \label{sec:cuspg29}

It follows from the analysis in~\cite{chl} that the lattice
$\CHL(G_{29},3)$ has a single conjugacy class of cusps. As mentioned
in the previous section, in order to see this, we need to consider
$G_{29}$-orbits of strata of the mirror arrangement given by
(irreducible) mirror intersections, and compute the number $\kappa_L$
(see p.~88 of~\cite{chl}). The cusps correspond to strata where
$\kappa_L=1$.

A list of representatives for all the $G$-orbits of irreducible mirror
intersections is given in Figure~\ref{tab:irred_g29}. Here, apart from
$L_{12343}$, we denote by $L_{j_1,\dots,j_r}$ the intersection of the
mirrors of $r_{j_1},\dots,r_{j_r}$, and we use the numbering given in
the diagram on page~\pageref{fig:cox}. For $L_{12343}$, we take the
intersection of the mirrors of $r_1$, $r_2$ and
$r_3r_4r_3=r_3r_4r_3^{-1}$.
\begin{table}
\begin{tabular}{|c||c|c|c|c|c|c|}
\hline
$L$        &      $L_1$      &          $L_{12}$            &     $L_{24}$       & $L_{123},L_{12343}$ &     $L_{124}$      & $L_{234}$\\
\hline
\# mirrors &       $1$       &            $3$               &        $4$         &       $6$           &        $9$         &  $12$\\
\hline
$\kappa_L$          & $1-\frac{2}{p}$ & $\frac{3}{2}(1-\frac{2}{p})$ & $2(1-\frac{2}{p})$ & $2(1-\frac{2}{p})$  & $3(1-\frac{2}{p})$ & $4(1-\frac{2}{p})$\\
\hline
\end{tabular}
\caption{Irreducible mirror intersections for $G_{29}$.}\label{tab:irred_g29}
\end{table}
For $p=3$, the only $G$-orbit with $\kappa_L=1$ is the $G$-orbit of
$L_{124}$. This means that the group $\Gamma=\CHL(G_{29},3)$ has a
single $\Gamma$-orbit of cusps, represented by the group generated by
$R_1,R_2$ and $R_4$. Note that this cusp group is a CHL group of the
form $\CHL(B_3,3)$, which fits in their analysis of parabolic cases.

We now work out the detailed structure of the cusp, using the general
framework of section~\ref{sec:cusps_general}.
First note that the reflections $R_1$, $R_2$ and $R_4$ fix a common
point in the ideal boundary $\partial_\infty H^3_{\mathbb{C}}$, given
in the basis used in section~\ref{sec:g29} by the vector
$$
%%v=(1,\zeta(\omega-1),0,-\zeta^2(1+i)).
v = (\zeta^2, \zeta^2 + 1, 0, \zeta^2 + \zeta  - 1).
$$ 
We use this as the first vector, and (a suitable multiple of)
$R_3v$ as the last basis vector. As the second vector, we use (a
suitable multiple of) the polar vector to the mirror of $R_1$, and as
the third one we take one that is orthogonal to both $v$ and $R_3v$
(and that makes the matrix of $R_2$ as simple as we could make it).

Concretely, we take
$$
Q = \left(\begin{matrix}
-\bar\omega   &  \bar\omega - 1 & i\sqrt{3}(1-\zeta)       &   i(3-2\sqrt{3})\\
1-\bar\omega  &     0           & i(\sqrt{3}-3)            & 3\zeta(\zeta-1)^2)\\
0             &     0           &            0             &  (\omega-1)(1+\omega\zeta)\\
\zeta+\omega  &     0           & \zeta(3-\sqrt{3})        &  i - \bar\omega - 5\zeta + 4
%  1                      & \zeta^2-2 & -\zeta^3+\zeta^2-\zeta+1 & -2\zeta^2+3\zeta-2\\
%  \zeta^3-2\zeta         &   0       & -\zeta^3+3\zeta^2-\zeta  & -3\zeta^2+6\zeta-3\\
%  0                      &   0       &            0             & -2\zeta^3+\zeta^2+\zeta-2\\
%  -\zeta^3-\zeta^2+\zeta &   0       & \zeta^3-3\zeta^2+\zeta   & -\zeta^3+4\zeta^2-4\zeta+1
\end{matrix}\right).
$$
Writing $S_j=Q^{-1}R_jQ$, we get
{\tiny
$$
S_1=\left(\begin{matrix}
  1 & 0      & 0 & 0\\
  0 & \omega & 1 & 0\\
  0 & 0      & 1 & 0\\
  0 & 0      & 0 & 1
\end{matrix}\right),
S_2=\left(\begin{matrix}
  1 & (\zeta-1)(1-\omega) & 3(1-\zeta)      & -2+(2-3\zeta)\omega\\
  0 & 1       & 0      & 0\\
  0 & -\omega & \omega & i-\omega\\
  0 & 0       & 0      & 1
\end{matrix}\right),
S_4=\left(\begin{matrix}
  1 & 0 & \zeta(\sqrt{3}-2)   & -2+\omega(2-3\zeta)\\
  0 & 1 &        -1           & \frac{1}{2}(i\sqrt{3}\zeta - \omega)\\
  0 & 0 &       \omega        & \frac{1}{2}i(3-\sqrt{3})\\
  0 & 0 &         0           & 1
\end{matrix}\right),
$$
}
which preserve the Hermitian form
\begin{equation}\label{eq:adaptedform-g29}
Q^*HQ=\left(\begin{matrix}
  0 &    0       &    0         &  1\\
  0 &    3       & \bar\omega-1 &  0\\
  0 & \omega-1   &    3         &  0\\
  1 &    0       &    0         &  0
\end{matrix}\right).
\end{equation}

In the notation of section~\ref{sec:cusps_general}, we have
$S_j=P(B_j,w_j,t_j)$ for 
\begin{equation}\label{eq:B-g29}
B_1=\left(\begin{matrix}
  \omega & 1\\
  0      & 1
\end{matrix}\right),\quad
B_2=\left(\begin{matrix}
  1       & 0\\
  -\omega & \omega
\end{matrix}\right),\quad
B_4=\left(\begin{matrix}
  1 & -1\\
  0 & \omega
\end{matrix}\right),
\end{equation}
and the translation parts are given by
\begin{equation}\label{eq:w-g29}
w_1=\left(\begin{matrix}0\\0 \end{matrix}\right),
w_2=\left(\begin{matrix}0\\i-\omega \end{matrix}\right),
w_4=\left(\begin{matrix}\frac{1}{2}(i\sqrt{3}\zeta - \omega)\\\frac{1}{2}i(3-\sqrt{3})\end{matrix}\right).
\end{equation}

One checks that the matrices $B_1$, $B_2$, $B_4$ generate a group
isomorphic to the Shephard-Todd group $G_5$ (which has order
72). Indeed, the matrices $B_2$ and $B_4$ are reflections of order 3
and $\br(B_2,B_4)=4$, so they generate a copy of $G_5$
(see~\cite{brmaro} for instance). One then checks that $B_1$ is in the
group generated by $B_2$ and $B_4$, for instance
$B_1=B_2[B_4^{-1},B_2]$.
%$B_1=B_2B_4^{-1}B_2B_4B_2^{-1}$.

This proves the following.
\begin{prop}\label{prop:linearpart-g29}
  The cusp of $\CHL(G_{29},3)$ is a central extension of a
  2-dimensional affine crystallographic group generated by
  reflections, with linear part $G_5$.
\end{prop}

We denote by $\mathcal{T}$ the subgroup of translations (i.e. the
unipotent subgroup) in the cusp. The first observation is that it
consists of matrices of the form $U(w,t)=P(Id,w,t)$ (see the notation
in section~\ref{sec:cusps_general}) such that $w$ has entries in
$Z[\omega]$.

This is not obvious, since the matrices $S_j$ have entries in
$\Z[\zeta]$. We start by enumerating some translations (that we
obtained by computing all words in $S_1,S_2,S_4$ of length at most 6).
The following can be checked by direct computation.
\begin{equation}\label{eq:tsl}
  \begin{array}{l}
   T_1=U\left((1,1),3(1-\frac{\sqrt{3}}{2})\right)=S_1S_2S_1^{-1}S_4^{-1}S_2^{-1}S_4\\
   T_2=U\left((1,0),-\frac{\sqrt{3}}{2}\right)=T_{(1,0)}=S_1S_2S_4^{-1}S_2^{-1}S_4S_2^{-1}\\
   T_3=U\left((\bar\omega,\bar\omega),-3(1-\frac{\sqrt{3}}{2})\right)=S_4^{-1}S_2^{-1}S_4S_1S_2S_1^{-1}\\
   T_4=U\left((0,\omega), -\frac{\sqrt{3}}{2}\right)=S_1S_4S_2S_1^{-1}S_4^{-1}S_2^{-1}
  \end{array}
\end{equation}
Note also that
$$
  [T_2,T_1]=U(0,\sqrt{3}),
$$
we denote the vertical translation by $V$.

Consider the abstract group with presentation
$$
G = \langle 
s_1,s_2,s_4
\, 
|
\,
s_1^3, s_2^3, s_4^3, \br_3(s_1,s_2), \br_4(s_2,s_4), [s_1,s_4], s_1s_4s_2s_1^{-1}s_4^{-1}s_2^{-1}
\rangle
$$ 
where the last two relations come from the right hand side of the
first and fourth row in equation~\eqref{eq:tsl}. We could have
included the other two relations coming from the second and third
rows, but it is easy to see that we would get the same abstract group
if we did.

One checks (most conveniently with a computer algebra system, say GAP)
that the abstract group $G$ has order 72, hence it is isomorphic to
$G_5$, which is the linear part of the complex crystallographic group
obtained by killing off the center of our cusp group.

As a consequence, we have the following.
\begin{prop}
  The translation subgroup $\mathcal{T}$ of the cusp group generated
  by $S_1,S_2,S_4$ is generated by the translations $T_1,T_2,T_3,T_4$
  in equation~\eqref{eq:tsl}.
\end{prop}
\begin{pf}
  Let $\Gamma_\infty$ denote the group generated by $S_1,S_2,S_4$ and
  let $\mathcal{T}_0$ denote the subgroup generated by
  $T_1,T_2,T_3,T_4$.
  
  The projection $\varphi:\Gamma_\infty\rightarrow U(K)$ given by
  mapping $P(B,w,t)$ to $B$ is a homomorphism onto $G_5$, and its
  kernel is the full translation subgroup $\mathcal{T}\subset
  \Gamma_\infty$.

  Every element in the kernel of $\varphi$ can be written as a product
  of conjugates of $S_1^3$, $S_1S_2S_1(S_2S_1S_2)^{-1}$,
  $(S_2S_4)^2(S_2S_4)^{-2}$, $S_1S_4S_1^{-1}S_4^{-1}$,
  $S_1S_4S_2S_1^{-1}S_4^{-1}S_2^{-1}$. Only the last element gives a
  non-trivial translation.

  Now one computes
  \begin{equation}
    \begin{array}{l}
      S_1T_4S_1^{-1} = (T_1T_3)^{-1} V\\
      S_2T_4S_2^{-1} = T_4^{-1}T_1^{-1}T_2 V^{-1}\\
      S_4T_4S_4^{-1} = T_2T_3 V^{-1},
    \end{array}
  \end{equation}
  which shows that the above conjugates are in fact in
  $\mathcal{T}_0$. In other words $\mathcal{T}=\mathcal{T}_0$.
\end{pf}

We now have the following.
\begin{prop}
  The complex crystallographic group obtained by projecting
  $\Gamma_\infty$ to $U(K)\rtimes \C^2$ has translation part
  $\Z[\omega]\times\Z[\omega]$. Its vertical translation group is
  generated by $V=U(0,\sqrt{3})$.
\end{prop}
\begin{pf}
  One easily checks that the vectors
  $u_1=(1,1),u_2=(1,0),u_3=(\bar\omega,\bar\omega),u_4=(0,\omega)$
  generate $\Z[\omega]\times\Z[\omega]$. 
\end{pf}

\subsection{Incommensurability} \label{sec:incomm}

We give an argument that was suggested by the referree, which replaces
our previous computational argument by a more geometric one.

\begin{prop}\label{prop:incomm}
  The lattices $\Gamma_1=\CHL(B_4,3,4)$ and $\Gamma_2\CHL(G_{29},3)$ are not commensurable.
\end{prop}
\begin{pf}
  Suppose they were commensurable, then there exists a $g\in GL(4,\C)$
  such that $g\Gamma_1 g^{-1}$ and $\Gamma_2$ have a common finite
  index subgroup $\Gamma$.  By irreducibility, this implies that
  $g^*H_1g=\lambda H_2$ for some $\lambda>0$, where $H_1$ and $H_2$
  are the Hermitian forms from equations~\eqref{eq:HcuspDM}
  and~\eqref{eq:adaptedform-g29}, respectively.

  We may assume that $\Gamma$ has a cusp represented by the standard
  basis vector $e_1$ and $H_2$ has the shape given in
  equation~\eqref{eq:adaptedform-g29}. The group $\Gamma$ contains
  both a nontrivial Heisenberg translation $gU_1(w_1,t_1)g^{-1}$ and a
  nontrivial vertical translation $gU_1(0,t_1')g^{-1}$, which can also
  be written as $U_2(w_2,t_2)$ and $U_2(0,t_2')$ (we use $U_j$ for the
  Heisenberg description using the Hermitian form $H_j$).
  
  Recall that the results of section~\ref{sec:cuspdm} imply that we
  may assume $w_1^*K_1w_1\in\Z$ and $t_1'\in i\Z$, and the results of
  section~\ref{sec:cuspg29} imply that we may assume (possibly after
  post-composing $g$ with some isometry of $H_2$), that
  $w_2^*K_2w_2\in\Z$ and $t_2'\in i\sqrt{3}\Z$.

  Moreover, the observation of equation~\eqref{eq:key} says that there
  is a $\lambda>0$ such that $w_2^*K_2w_2=\lambda w_1^*K_1w_1$ and
  $t_2'=\lambda t_1'$.

  The first equation implies $\lambda\in\Q$, whereas the second
  implies $\lambda\notin\Q$, contradiction.
\end{pf}

\appendix
\section{Coxeter diagrams}\label{sec:cox}

For the reader's convenience, we gather Coxeter diagrams for the
Shephard-Todd groups. These were worked out by
Coxeter~\cite{coxeterbook} and Shephard-Todd~\cite{shephardtodd}, see
also the Appendix 2 in~\cite{brmaro} for a convenient list.

\begin{figure}
  \hfill
  \subfigure[$A_4$]{
    \includegraphics[height=0.5cm]{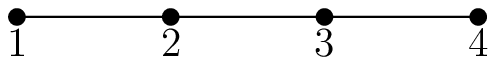}
  }
  \hfill
  \subfigure[$B_4$]{
    \includegraphics[height=0.5cm]{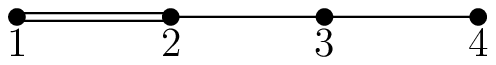}
  }
  \hfill\,
  \\
  \hfill
  \subfigure[$G_{28}$]{
    \includegraphics[height=0.5cm]{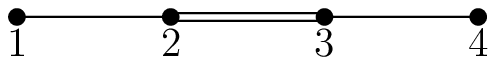}
  }
  \hfill
  \subfigure[$G_{29}$]{
    \includegraphics[height=2cm]{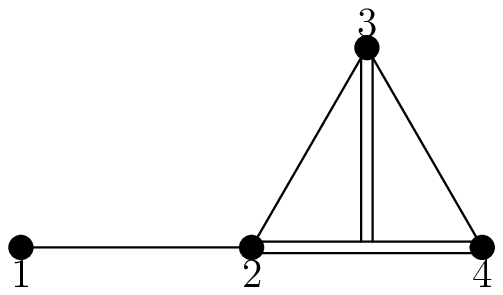}
  }
  \hfill\,
  \\
  \hfill
  \subfigure[$G_{30}$]{
    \includegraphics[height=0.7cm]{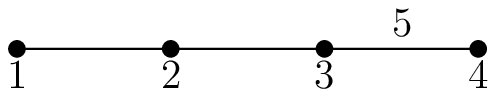}
  }
  \hfill
  \subfigure[$G_{31}$]{
    \includegraphics[height=1.8cm]{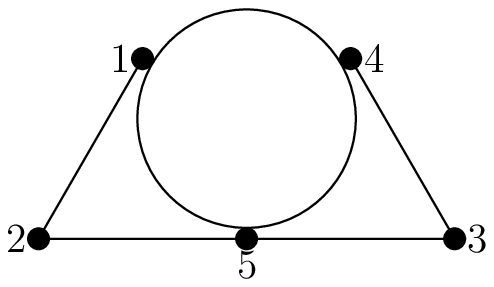}
  }
  \hfill\,
  \\
  \hfill
  \subfigure[$G_{33}$]{
    \includegraphics[height=1.8cm]{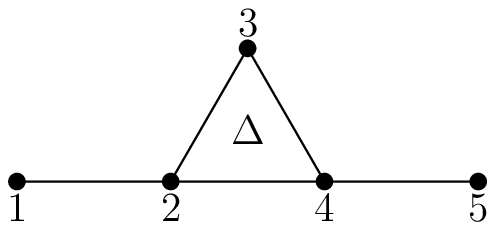}
  }
  \hfill
  \subfigure[$G_{34}$]{
    \includegraphics[height=1.8cm]{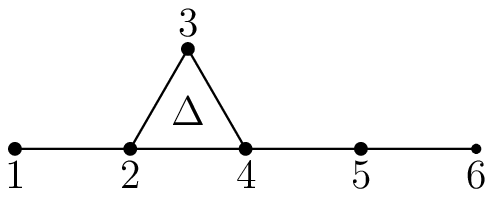}
  }
  \hfill\,
  \\
  \hfill
  \subfigure[$G_{35}$]{
    \includegraphics[height=1.8cm]{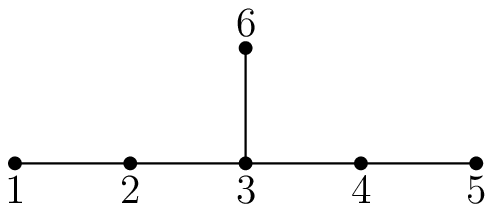}
  }
  \hfill
  \subfigure[$G_{36}$]{
    \includegraphics[height=1.8cm]{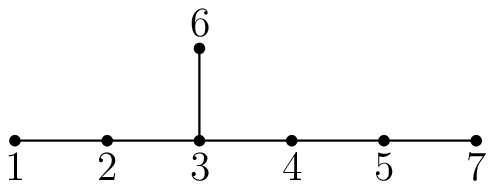}
  }
  \hfill\,
  \\
  \hfill
  \subfigure[$G_{37}$]{
    \includegraphics[height=1.8cm]{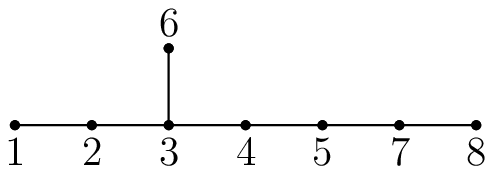}
  }\hfill\,
\caption{Coxeter diagrams for exceptional Shephard-Todd groups. The
  symbol $\Delta$ in the diagrams for $G_{33}$ and $G_{34}$ stand for
  the relation $(R_2R_4R_3)^2=(R_4R_3R_2)^2=(R_3R_2R_4)^2$.}\label{fig:cox}
\end{figure}

\section{Rough commensurability invariants}\label{sec:rough}

In this section, we gather rough commensurability invariants for the
Couwenberg-Heckman-Looijenga lattices in $PU(n,1)$ that correspond to
the exceptional Shephard-Todd groups. For each such Shephard-Todd
group, we list values of the order of reflections that yield lattices,
and mention whether the corresponding lattices are cocompact (C/NC)
and arithmetic (A/NA). We also give their adjoint trace field $\Q(\tr
Ad \Gamma)$.

\begin{table}
\begin{tabular}{|c|c|c|}
\hline
Shephard-Todd & Other description & $p$\\
\hline
$G_{23}$      & Coxeter $H_3$, $\Sc(\sigma_{10},p)$ & ${\bf 3},{\bf 4},{\bf 5},{\bf 10}$\\
\hline
$G_{24}$      &  $\Sc(\overline{\sigma}_4,p)$       & ${\bf 3},{\color{red} 4},{\color{red}\bf 5},{\color{red}6},{\color{red}\bf 8},{\color{red}\bf 12}$\\
\hline
$G_{27}$      & $\T({\bf S_2},p)$                   & ${\bf 3}, {\color{red}4}, {\color{red}\bf 5}$\\
\hline
\end{tabular}
\caption{2-dimensional CHL lattices. Note that $G_{25}$ and $G_{26}$
  are not listed because they yield Deligne-Mostow groups. Values of
  $p$ that appear in bold-face correspond to cocompact lattices, the
  ones that appear in red give non-arithmetic lattices.} \label{tab:2d}
\end{table}

\begin{table}
\begin{tabular}{|c|c|c|c|c|c|c|}
\hline
Shephard-Todd & Dimension & Other description                   & $p$ or $(p_1,p_2)$ & A? & C? & $\Q(\tr Ad\Gamma)$\\
\hline
$G_{28}$      &     3     &    Coxeter $F_4$  & ${(2,4)}$                  & A  & NC & $\Q$\\
              &           &                   & ${\bf (2,5)}$              & A  & C  & $\Q(\sqrt{5})$\\
              &           &                   & ${(2,6)}$                  & A  & NC & $\Q$\\
              &           &                   & ${\bf (2,8)}$              & A  & C  & $\Q(\sqrt{2})$\\
              &           &                   & ${\bf (2,12)}$             & A  & C  & $\Q(\sqrt{3})$\\
              &           &                   & ${(3,3)}$                  & A  & NC & $\Q$\\
              &           &                   & ${\bf (3,4)}$              & A  & C  & $\Q(\sqrt{3})$\\
              &           &                   & ${(3,6)}$                  & A  & NC & $\Q$\\
              &           &                   & ${\bf (3,12)}$             & A  & C  & $\Q(\sqrt{3})$\\
              &           &                   & ${(4,4)}$                  & A  & NC & $\Q$\\
              &           &                   & ${(6,6)}$                  & A  & NC & $\Q$\\
\hline
$G_{29}$      &    3      &                   & ${\color{red}3}$           & NA & NC & $\Q(\sqrt{3})$\\
              &           &                   & ${4}$                      & A  & NC & $\Q$\\
\hline
$G_{30}$      &    3      & Coxeter $H_4$     & ${\bf 3}$                  & A  & C  & $\Q(\sqrt{5})$\\
              &           &                   & ${\bf 5}$                  & A  & C  & $\Q(\sqrt{5})$\\
\hline
$G_{31}$      &    3      &                   & ${3}$                      & A  & NC & $\Q$\\
              &           &                   & ${\bf 5}$                  & A  & C  & $\Q(\sqrt{5})$\\
\hline
$G_{33}$      &    4      &                   & $3$                        & A  & NC & $\Q$\\
\hline
$G_{34}$      &    5      &                   & $3$                        & A  & NC & $\Q$\\
\hline
$G_{35}$      &    5      &    $W(E_6)$       & $3$                        & A  & NC & $\Q$\\
              &           &                   & $4$                        & A  & NC & $\Q$\\
\hline
$G_{36}$      &    6      &    $W(E_7)$       & $3$                        & A  & NC & $\Q$\\
\hline
$G_{37}$      &    7      & $W(E_8)$          & $3$                        & A  & NC & $\Q$       \\
\hline
\end{tabular}
\caption{Rough commensurability invariants for CHL lattices. The group
  $G_{32}$ is not listed, since it is also $A_3$ and yields
  Deligne-Mostow groups.}\label{tab:list}
\end{table}

\end{document}